\input ./style/arxiv-general.cfg
\documentclass[aos,MSNbibl,nameyear,seceqn,dvips]{arximspdf}
\makeatletter
   \@ifpackageloaded{graphicx}{}{\usepackage{graphicx}}
\makeatother
\usepackage{dcolumn}

\doi{10.1214/15-AOS1367}
\volume{44}
\issue{1}
\pubyear{2016}
\firstpage{288}
\lastpage{317}
\docsubty{FLA}

\makeatletter
\newcolumntype{d}[1]{D{.}{.}{#1}}
\renewcommand{\mid}{|}
\newcommand{\rrvert}{\vert}
\newcommand{\rrVert}{\Vert}
\newcommand{\llvert}{\vert}
\newcommand{\llVert}{\Vert}
\renewcommand{\mid}{|}
\newcommand{\mathds}{\mathbb}
\def\convInD{\stackrel{d}{\rightarrow}}
\newcommand{\vbeta}{\bolds{\beta}}
\newcommand{\vc}{\mathbf{c}}
\newcommand{\vW}{\mathbf{W}}
\newcommand{\vt}{\mathbf{t}}
\newcommand{\vs}{\mathbf{s}}
\newcommand{\vx}{\mathbf{x}}
\newcommand{\vz}{\mathbf{z}}
\newcommand{\vgamma}{\bolds{\gamma}}
\newcommand{\vxi}{\bolds{\xi}}
\newcommand{\vpi}{\bolds{\pi}}
\newcommand{\vtheta}{\bolds{\theta}}
\newcommand{\vdelta}{\bolds{\delta}}
\newcommand{\veta}{\bolds{\eta}}
\newcommand{\vkappa}{\bolds{\kappa}}
\newcommand{\vnull}{{\mathbf0}}
\newtheorem{theorem}{Theorem}[section]
\newtheorem{lemma}{Lemma}
\newtheorem{corollary}{Corollary}
\newproclaim{condition}{Condition}
\newproclaim{defi}{Definition}
\newproclaim{rem}{Remark}
\makeatother

\begin{document}
\begin{frontmatter}

\title{Partially linear additive quantile regression in~ultra-high dimension}
\runtitle{Ultra-high dimensional PLA quantile regression}

\begin{aug}
\author[A]{\fnms{Ben}~\snm{Sherwood}\corref{}\ead[label=e1]{bsherwo2@jhu.edu}}
\and
\author[B]{\fnms{Lan}~\snm{Wang}\thanksref{T2}\ead[label=e2]{wangx346@umn.edu}}
\runauthor{B. Sherwood and L. Wang}
\affiliation{Johns Hopkins University and University of Minnesota}
\address[A]{Department of Biostatistics\\
Johns Hopkins University\\
Baltimore, Maryland 21205\\
USA\\
\printead{e1}}
\address[B]{School of Statistics\\
University of Minnesota\\
Minneapolis, Minnesota 55455\\
USA\\
\printead{e2}}
\end{aug}
\thankstext{T2}{Supported in part by NSF Grant DMS-13-08960.}

%
\received{\smonth{9} \syear{2014}}
%
\revised{\smonth{7} \syear{2015}}

%
\begin{abstract}
We consider a flexible semiparametric quantile regression model for analyzing
high dimensional heterogeneous data. This model has
several appealing features: (1) By considering different
conditional quantiles, we may \mbox{obtain} a more complete picture of the
conditional distribution of a response variable given
high dimensional covariates.
(2) The sparsity level is allowed to be different at different quantile levels.
(3) The partially linear additive structure accommodates nonlinearity
and circumvents the curse of dimensionality. (4) It is naturally robust to
heavy-tailed distributions. In this paper, we approximate the
nonlinear components using B-spline basis functions. We first study
estimation under this model when the nonzero components
are known in advance and the number of covariates in the linear part diverges.
We then investigate a nonconvex penalized estimator for simultaneous
variable selection and estimation. We derive its oracle property for
a general class of nonconvex penalty functions
in the presence of ultra-high dimensional covariates under relaxed conditions.
To tackle the challenges of nonsmooth loss function, nonconvex penalty function
and the presence of nonlinear components, we combine a recently developed
convex-differencing method with 
modern empirical process techniques. 
Monte Carlo simulations and an application to a microarray study
demonstrate the effectiveness of the proposed method. We also discuss
how the method
for a single quantile of interest can be extended to simultaneous
variable selection and estimation at multiple quantiles.
\end{abstract}

%
\begin{keyword}[class=AMS]
\kwd[Primary ]{62G35}
\kwd[; secondary ]{62G20}
\end{keyword}
\begin{keyword}
\kwd{Quantile regression}
\kwd{high dimensional data}
\kwd{nonconvex penalty}
\kwd{partial linear}
\kwd{variable selection}
\end{keyword}
\end{frontmatter}

\section{Introduction}\label{sec1}
In this article, we study a flexible partially linear additive
quantile regression model for analyzing high dimensional data.
For the $i$th subject, we observe $ \{Y_i,\vx_i,\vz_i \}$,
where $\vx_i= (x_{i1},\ldots,x_{ip_n})'$ is a $p_n$-dimensional
vector of covariates
and $\vz_i= (z_{i1},\ldots,z_{id})'$ is a $d$-dimensional vector of
covariates, $i=1,\ldots,n$.
The $\tau$th ($0<\tau<1$) conditional quantile of $Y_i$ given $\vx_i$,
$\vz_i$ is defined as $Q_{Y_i|\vx_i,\vz_i}(\tau)=\inf\{t: F(t|\vx_i
,\vz_i)\geq\tau\}$,
where $F(\cdot|\vx_i,\vz_i)$ is the conditional distribution function
of $Y_i$ given $\vx_i$
and $\vz_i$. The case $\tau=1/2$ corresponds to the conditional median.
We consider the following semiparametric model for the conditional
quantile function
%
%
\begin{equation}
\label{model} Q_{Y_i|\vx_i,\vz_i}(\tau)=\vx_i'
\vbeta_0+g_0(\vz_i),
\end{equation}
where $g_0(\vz_i)=g_{00}+\sum_{j=1}^dg_{0j}(z_{ij})$, with $g_{00}\in
\mathcal{R}$.
It is assumed that $g_{0j}$ satisfy $E(g_{0j}(z_{ij}))=0$ for
identification purposes.
Let $\varepsilon_i=Y_i-Q_{Y_i|\vx_i,\vz_i}(\tau)$, then $\varepsilon_i$
satisfies $P(\varepsilon_i\leq0|\vx_i,\vz_i)=\tau$ and we may also write
$Y_i=\vx_i'\vbeta_0+g_0(\vz_i)+\varepsilon_i$.
In the rest of the paper, we will drop the dependence on $\tau$
in the notation for simplicity.

Modeling conditional quantiles in high dimension is of significant
importance for several reasons.
First, it is well recognized that high dimensional data are often heterogeneous.
How the covariate influence the center of the conditional distribution
can be very different from how they influence the tails.
As a result, focusing on the conditional mean function alone can be misleading.
By estimating conditional quantiles at different quantile levels,
we are able to gain a more complete picture of the relationship between
the covariates and the response variable.
Second, in the high dimensional setting,
the quantile regression framework also allows a more realistic
interpretation of the sparsity
of the covariate effects, which we refer to as quantile-adaptive sparsity.
That is, we assume a small subset of covariates influence the conditional
distribution. However, when we estimate different conditional
quantiles, we allow the
subsets of active covariates to be different [\citet{quantUltraHigh};
\citet{heWang}].
Furthermore, the conditional quantiles are often of direct interest to
the researchers.
For example, for the birth weight data we analyzed in Section~\ref{sec5}, low
birth weight, which corresponds to the low tail of the conditional
distribution, is of direct interest to the doctors. Another advantage
of quantile regression is that it is naturally robust to outlier
contamination associated with heavy-tailed errors. For high dimensional
data, identifying outliers can be difficult. The robustness of quantile
regression provides a certain degree of protection.


Linear quantile regression with high dimensional covariates was
investigated by
\citeauthor{B1} [(\citeyear{B1}), Lasso penalty] and
Wang, Wu and Li [(\citeyear{quantUltraHigh}), nonconvex penalty].
The partially linear additive structure we consider in this paper is
useful for incorporating nonlinearity
in the model while circumventing the curse of dimensionality.
We are interested in the case
$p_n$ is of a similar order of $n$ or much larger than $n$.
For applications in microarray data analysis,
the vector $\vx_i$ often contains the measurements on thousands of
genes, while the
vector $\vz_i$ contains the measurements of clinical or environment
variables, such as
age and weight. For example, in the birth weight example of Section~\ref{sec5},
mother's age is modeled nonparametrically
as exploratory analysis reveals a possible nonlinear effect. In
general, model specification
can be challenging in high dimension; see Section~\ref{sec7} for some
further discussion.

We approximate the nonparametric components using B-spline basis
functions, which are computationally convenient
and often accurate. First, we study the asymptotic theory
of estimating the model (\ref{model}) when $p_n$ diverges.
In our setting, this corresponds to the oracle model, that is, the
one we obtain if we know which covariates are important in advance.
This is along the line of the work
of \citet{Welsh}, \citet{BaiWu} and \citet{HeShao}
for $M$-regression with diverging number of parameters and
possibly nonsmooth objective functions, which, however,
were restricted to linear regression. \citet{partMeanVaryCoef}
derived the asymptotic theory of profile kernel estimator
for general semiparametric models with diverging number
of parameter while assuming a smooth quasi-likelihood function.
Second, we propose a nonconvex
penalized regression estimator when $p_n$ is of an exponential order of
$n$ and
the model has a sparse structure. For a general class of nonsmooth
penalty functions, including the popular SCAD [\citet{fanLi}]
and MCP [\citet{Zhang}] penalty, we derive the oracle property of the
proposed estimator under
relaxed conditions. An interesting finding is that
solving the nonconvex penalized estimator can be achieved via
solving a series of weighted quantile regression problems, which can
be conveniently implemented using existing software packages.

Deriving the asymptotic properties of the penalized estimator is very
challenging as
we need to simultaneously deal with the nonsmooth loss function, the
nonconvex penalty
function, approximation of nonlinear functions and very high \mbox{dimensionality}.
To tackle these challenges, we combine a recently developed
convex-differencing method with 
modern empirical process techniques. The method relies on a
representation of the penalized loss function as the difference
of two convex functions, which leads to a sufficient local
optimality condition [\citet{TaoAn97}, \citet{quantUltraHigh}].
Empirical process techniques are
introduced to derive various error bounds associated with
the nonsmooth objective function which contains both high dimensional
linear covariates and approximations of nonlinear components.
It is worth pointing out that our approach is different from what was
used in the
recent literature for studying the theory of
high dimensional semiparametric mean regression and is able to
considerably weaken
the conditions required in the literature. In particular, we do not
need moment
conditions for the random error and allow it to depend on the covariates.

Existing work on penalized semiparametric regression has been largely
limited to mean regression with
fixed $p$;
see, for example, \citet{Bunea}, \citet{LiangLi}, \citet{WX},
\citet{scadPartLinMean}, \citet{KLZ} and
\citet{WLLC}.
Important progress in the high dimensional $p$ setting
has been recently made by
\citeauthor{annalsMeanParLin} [(\citeyear{annalsMeanParLin}), still
assumes $p<n$] for
partially linear regression, \citet{HHW} for
additive models,
\citeauthor{Li} [(\citeyear{Li}), $p=o(n)$] for semivarying coefficient models,
among others. When $p$ is fixed, the semiparametric quantile regression
model was considered by
\citet{heShi96}, \citet{heQuantLong}, \citet{quantVaryModel},
among others. \citet{Tang} considered a two-step procedure for a
nonparametric varying
coefficients quantile regression model with a diverging number of
nonparametric functional coefficients.
They required two separate tuning parameters and quite complex design
conditions.

The rest of this article is organized as follows. In Section~\ref
{sec2}, we
present the partially linear additive quantile regression model and
discuss the properties of the oracle estimator.
In Section~\ref{sec3}, we present a nonconvex penalized method for simultaneous
variable selection and estimation and derive its
oracle property. In Section~\ref{sec4}, we assess the performance of the
proposed penalized
estimator via Monte Carlo simulations.
We analyze a birth weight data set while accounting for gene expression
measurements in Section~\ref{sec5}.
In Section~\ref{sec6}, we consider an extension to simultaneous
estimation and
variable selection at multiple quantiles.
Section~\ref{sec7} concludes the paper with a discussion of related issues.
The proofs are given in the \hyperref[appe]{Appendix}.
Some of the technical details and additional numerical results are
provided in online supplementary material [\citet{Supp}].

\section{Partially linear additive quantile regression with diverging
number of parameters}\label{sec2}
For high dimensional inference, it is often assumed that the vector of
coefficients
$\vbeta_0=(\beta_{01},\beta_{02},\ldots,\beta_{0p_n})'$ 
in model (\ref{model}) is sparse, that is, most of its components are
zero. Let $A = \{1 \leq j \leq p_n: \beta_{0j} \neq0\}$ be the index
set of nonzero coefficients and $q_n = |A|$ be the cardinality of $A$.
The set $A$ is unknown and will be estimated.
Without loss of generality, we assume that the first $q_n$ components of
$\vbeta_0$ are nonzero and the remaining $p_n-q_n$ components are
zero. 
Hence, we can write $\vbeta_0 = (\vbeta_{01}',\mathbf
{0}_{p_n-q_n}' )'$, where $\mathbf{0}_{p_n-q_n}$ denotes the
$(p_n-q_n)$-vector of zeros.
Let $X$ be the $n \times p_n$ matrix of linear covariates and write it
as $X = (X_1,\ldots, X_{p_n})$.
Let $X_A$ be the submatrix consisting of the first $q_n$ columns of $X$
corresponding to the active covariates.
For technical simplicity, we assume $x_{i}$ is centered to have mean zero;
and $z_{ij} \in[0,1]$, $\forall i, j$.

\subsection{Oracle estimator}\label{sec2.1}
We first study the estimator we would obtain when the index set $A$ is
known in advance, which we
refer to as the oracle estimator. Our asymptotic framework allows
$q_n$, the size of $A$, to increase with $n$.
This resonates with the perspective that a more complex statistical
model can be fit when more data are collected.

We use a linear combination of B-spline basis functions to approximate
the unknown nonlinear functions
$g_0(\cdot)$.
To introduce the B-spline functions, we start with two definitions.

\begin{defi*}
Let $r
\equiv m + v$,
where $m$ is a positive integer and $v\in(0,1]$. Define $\mathcal
{H}_r$ as the collection of functions $h(\cdot)$ on $[0,1]$ whose
$m$th derivative $h^{(m)}(\cdot)$ satisfies the H\"{o}lder condition
of order $v$. That is, for any $h(\cdot) \in\mathcal{H}_r$, there
exists some positive constant $C$
such that
%
%
\begin{equation}
\label{holder_cond} \bigl\llvert h^{(m)}\bigl(z'\bigr) -
h^{(m)}(z)\bigr\rrvert\leq C \bigl\llvert z' - z\bigr
\rrvert^v\qquad\forall0 \leq z', z \leq1.
\end{equation}
\end{defi*}

Assume for some $r \geq1.5$, the nonparametric
component $g_{0k}(\cdot) \in\mathcal{H}_r$.
Let $\vpi(t)= (b_1(t),\ldots,b_{k_n+l+1}(t) )'$ denote
a vector of normalized B-spline basis functions of order $l+1$
with $k_n$ quasi-uniform internal knots on $[0,1]$.
Then $g_{0k}(\cdot)$ can be approximated using a linear combination of
B-spline basis functions in
$\bolds{\Pi}(\vz_i)=(1,\vpi(z_{i1})',\ldots,\vpi(z_{id})')'$.
We refer to \citet{Schumaker} for details of the B-spline construction,
and the result
that there exists $\vxi_{0}\in\mathcal{R}^{L_n}$, where
$L_n=d(k_n+l+1)+1$, such that
$\sup_{\vz_i}|\bolds{\Pi}(\vz_i)'\vxi_{0} - g_{0}(\vz_i)| =
O
(k_n^{-r} )$.
For ease of notation and simplicity of proofs, we use the same number
of basis functions for
all nonlinear components in model (\ref{model}). In practice, such
restrictions are not necessary.

Now consider quantile regression with the oracle information that the
last $(p_n-q_n)$
elements of $\vbeta_0$ are all zero. Let
%
%
\begin{equation}
\label{orcObjFun} (\hat{\vbeta}_1, \hat{\vxi} ) = \mathop{\operatorname
{argmin}}_{ (\vbeta
_1, \vxi)} \frac{1}{n} \sum
_{i=1}^n\rho_\tau\bigl(Y_i
- \vx_{A_i} '\vbeta_1 - \bolds{\Pi}(
\vz_i)'\vxi\bigr),
\end{equation}
%
where\vspace*{-2pt}
$\rho_\tau(u) = u(\tau- I(u<0))$ is the quantile loss function and
$\vx_{A_1}',\ldots,\vx_{A_n}'$ denote the row vectors of $X_A$.
The\vspace*{-2pt} oracle estimator for $\vbeta_0$ is $ (\hat{\vbeta}_1',
\mathbf{0}_{p_n-q_n}' )'$.
Write $\hat{\vxi}=(\hat{\xi}_0,\hat{\vxi}_1',\ldots,\hat{\vxi}_d')'$
where $\hat{\xi}_0\in\mathcal{R}$ and $\hat{\vxi}_j\in\mathcal
{R}^{k_n+l+1}$, $j=1\ldots,d$.
The estimator
for the nonparametric function $g_{0j}$ is
\[
\hat{g}_j(z_{ij}) = \vpi(z_{ij})'
\hat{\vxi}_j-n^{-1}\sum_{i=1}^n
\vpi(z_{ij})'\hat{\vxi}_j,
\]
for $j=1,\ldots,d$; for $g_{00}$
is $\hat{g}_{0}=\hat{\xi}_0+n^{-1}\sum_{i=1}^n\sum_{j=1}^d\vpi(z_{ij})
'\hat{\vxi}_j$.
The centering of $\hat{g}_{j}$ is the sample analog of the identifiability
condition $E[g_{0j}(\vz_i)] = 0$. The estimator of $g_{0}(\vz_i)$ is
$\hat{g}(\vz_i)=\hat{g}_{0}+\sum_{j=1}^d\hat{g}_{j}(z_{ij})$.

\subsection{Asymptotic properties}\label{sec2.2}
We next present the asymptotic properties of the oracle estimators as
$q_n$ diverges.

\begin{defi*}
Given $\vz
= (z_1,\ldots
,z_d)'$, the function $g(\vz)$ is said to belong to the class of
functions $\mathcal{G}$ if it has the representation $g(\vz) = \alpha
+\sum_{k=1}^d g_k(\vz_k)$, $\alpha\in\mathcal{R}$,
$g_k \in\mathcal{H}_r$ and $E[g_k(\vz_{k})]=0$.
\end{defi*}

Let
\[
h^*_j(\cdot) = \mathop{\operatorname{arg}\operatorname
{inf}}_{h_j(\cdot) \in\mathcal{G}} \sum_{i=1}^nE \bigl[
f_i(0) \bigl(x_{ij} - h_j(\vz_i)
\bigr)^2 \bigr],
\]
where $f_i(\cdot)$ is the probability density function of $\varepsilon
_i$ given $(\vx_i,\vz_i)$.
Let $m_j(\vz) = E [x_{ij} \mid\vz_i=\vz]$, then it can\vspace*{1pt}
be shown that
$h^*_j(\cdot)$ is the weighted projection of
$m_j(\cdot)$ into $\mathcal{G}$ under the $L_2$ norm, where the
weights $f_i(0)$ are included to account for the possibly heterogeneous errors.
Furthermore, let $x_{A_{ij}}$ be the $(i,j)$th element of $X_A$. Define
$\delta_{ij} \equiv x_{A_{ij}} - h^*_j(\vz_i)$, $\vdelta_{i} =
(\delta_{i1},\ldots,\delta_{iq_n} )' \in\mathcal{R}^{q_n}$
and $\Delta_n = (\vdelta_1,\ldots,\vdelta_n )' \in
\mathds{R}^{n \times q_n}$. Let $H$ be the $n \times q_n$ matrix with
the $(i,j)$th element $H_{ij}= h_j^*(\vz_i)$, then $X_A = H + \Delta
_n$.

The following technical conditions are imposed for analyzing the
asymptotic behavior of $\hat{\vbeta}_1$ and $\hat{g}$.

%
\begin{condition}[(Conditions on the random error)]\label{cond_f}
The random error $\varepsilon_i$ has the conditional distribution
function $F_i$ and continuous
conditional density function $f_i$, given $\vx_i$, $\vz_i$.
The $f_i$ are uniformly bounded away from 0 and
infinity in a neighborhood of zero, its first derivative $f_i'$ has a
uniform upper bound in a neighborhood of zero,
for $1\leq i\leq n$.
\end{condition}

%
\begin{condition}[(Conditions on the covariates)]\label{highd_cond_x}
There exist positive constants $M_1$ and $M_2$ such that $|x_{ij}| \leq
M_1$, $\forall1\leq i \leq n, 1 \leq j \leq p_n$ and
$E[\delta_{ij}^{4}] \leq M_2$, $\forall1\leq i \leq n, 1 \leq j
\leq q_n$.
There exist finite positive constants $C_1$ and $C_2$ such that with
probability one
%
\[
C_1 \leq\lambda_{\max} \bigl( n^{-1}X_AX_A'
\bigr) \leq C_2, \qquad C_1 \leq\lambda_{\max}
\bigl(n^{-1}\Delta_n \Delta_n'
\bigr) \leq C_2.
\]
\end{condition}

%
\begin{condition}[(Condition on the nonlinear functions)]\label{cond_g_h}
For $r=m+v> 1.5$ $g_0 \in\mathcal{G}$.
\end{condition}


%
\begin{condition}[(Condition on the B-spline basis)]
\label{cond_j_n}
The\vspace*{1pt} dimension of the spline basis $k_n$ has the following rate
$k_n \approx n^{1/(2r+1)}$.
\end{condition}

%
\begin{condition}[(Condition on model size)]\label{cond_sigma_large_p}
$q_n = O (n^{C_3} )$ for some $C_3< \frac{1}{3}$.%
\end{condition}



Condition \ref{cond_f} is considerably more relaxed than what is
usually imposed on the random error for
the theory of high dimensional mean regression, which often requires
Gaussian or sub-Gaussian tail condition.
Condition~\ref{highd_cond_x} is about the behavior of the covariates and the design
matrix under the oracle model, which
is not restrictive.
Condition \ref{cond_g_h} is typical for the application of B-splines.
\citet{Stone85} showed that B-splines basis
functions can be used to effectively approximate functions satisfying
H\"{o}lder's condition. 
Condition \ref{cond_j_n} provides the rate of $k_n$ needed for the
optimal convergence rate of $\hat{g}$.
Condition~\ref{cond_sigma_large_p} is standard for linear models with diverging number of parameters.\vadjust{\goodbreak}

The following theorem summarizes the asymptotic properties of the
oracle estimators.

%
\begin{theorem}
\label{large_q_oracle}
Assume Conditions \ref{cond_f}--\ref{cond_sigma_large_p} hold.
Then
\begin{eqnarray*}
\llVert\hat{\bolds{\vbeta}}_1-\vbeta_{01}\rrVert&=&
O_p \bigl(\sqrt{n^{-1}q_n} \bigr),
\\
n^{-1} \sum_{i=1}^n \bigl(
\hat{g}(\vz_i) - g_0(\vz_i)
\bigr)^2 &=& O_p \bigl(n^{-1}(q_n+k_n)
\bigr).
\end{eqnarray*}
\end{theorem}

An interesting observation is that since we allow $q_n$ to diverge with
$n$, it influences the rates for estimating both $\vbeta$ and $g$. As
$q_n$ diverges, to investigate the asymptotic distribution of $\hat
{\vbeta}_1$, we consider estimating an arbitrary linear combination of
the components
of $\vbeta_{01}$.

%
\begin{theorem}
\label{large_q_clt}
Assume the conditions of Theorem \ref{large_q_oracle} hold.
Let $A_n$ be an $l\times q_n$ matrix with $l$ fixed and
$A_nA_n'\rightarrow G$, a positive definite matrix, then
\[
\sqrt{n}A_n\Sigma_n^{-1/2} (\hat{
\vbeta}_1-\vbeta_{01} ) \rightarrow N(\mathbf{0}_l,G)
\]
in distribution, where $\Sigma_n=K_n^{-1}S_n K_n^{-1}$ with
$K_n=n^{-1}\Delta_n'B_n \Delta_n$,
$S_n =\break n^{-1}\tau(1-\tau) \Delta_n'\Delta_n$,
and $B_n=\operatorname{diag}(f_1(0),\ldots,f_n(0))$ is an $n \times n$
diagonal matrix with
$f_i(0)$ denoting the conditional density function of $\varepsilon_i$
given $(\vx_i,\vz_i)$ evaluated at zero.
\end{theorem}

If we consider the case where $q$ is fixed and finite, then we have the
following result
regarding the behavior of the oracle estimator.

%
\begin{corollary}
\label{fixed_q_clt}
Assume $q$ is a fixed positive integer,
$n^{-1}\Delta_n'B_n \Delta_n\rightarrow\Sigma_1$ and
$n^{-1}\tau(1-\tau)\Delta_n'\Delta_n\rightarrow\Sigma_2$, where
$\Sigma_1$ and $\Sigma_2$
are positive definite matrices.
If Conditions \ref{cond_f}--\ref{cond_j_n} hold, then
\begin{eqnarray*}
\sqrt{n} (\hat{\vbeta}_1-\vbeta_{01} ) &\convInD& N
\bigl(\mathbf{0}_q, \Sigma_1^{-1}
\Sigma_2\Sigma_1^{-1} \bigr),
\\
n^{-1}\sum_{i=1}^n \bigl(
\hat{g}(\vz_i) - g_0(\vz_i)
\bigr)^2 &=& O_p \bigl(n^{-2r/(2r+1)} \bigr).
\end{eqnarray*}
\end{corollary}

In the case $q_n$ is fixed, the rates reduce to the classical
$n^{-1/2}$ rate for estimating $\vbeta$ and $n^{-2r/(2r+1)}$ for
estimating $g$, the latter
which is consistent with \citet{Stone85} for the optimal rate of
convergence.

\section{Nonconvex penalized estimation for partially linear additive
quantile regression with ultra-high dimensional covariates}\label{sec3}

\subsection{Nonconvex penalized estimator}\label{sec3.1}
In real data analysis, we do not know which of the $p_n$ covariates in
$\vx_i$ are important.
To encourage sparse\vadjust{\goodbreak} estimation, we minimize the following penalized
objective function for estimating $(\vbeta_0, \vxi_0)$,
%
%
\begin{equation}
\label{penObj} Q^P(\vbeta,\vxi) = n^{-1} \sum
_{i=1}^n\rho_\tau\bigl(Y_i -
\vx_i'\vbeta- \bolds{\Pi}(\vz_i)'
\vxi\bigr) + \sum_{j=1}^{p_n}
p_\lambda\bigl(\llvert\beta_j\rrvert\bigr),
\end{equation}
where $p_\lambda(\cdot)$ is a penalty function with tuning parameter
$\lambda$.
The $L_1$ penalty or Lasso [\citet{T1}] is a popular choice for
penalized estimation. However,
the $L_1$ penalty is known to over-penalize large coefficients, tends
to be biased
and requires strong conditions on the design matrix to achieve
selection consistency. 
This is usually not 
a concern for prediction, but can be undesirable if the goal is to
identify the underlying model. In comparison, an appropriate nonconvex
penalty function can effectively overcome this problem [\citet
{fanLi}]. In this paper, we consider
two such popular choices of penalty functions: the SCAD [\citet
{fanLi}] and MCP [\citet{Zhang}] penalty functions. For the SCAD
penalty function,
\begin{eqnarray*}
p_\lambda\bigl(\llvert\beta\rrvert\bigr) &=& \lambda\llvert\beta
\rrvert I\bigl(0 \leq\llvert\beta\rrvert< \lambda\bigr) + \frac{a
\lambda\llvert \beta\rrvert - (\beta^2 + \lambda^2
)/2}{a-1}I
\bigl(\lambda\leq\llvert\beta\rrvert\leq a\lambda\bigr)
\\
&&{}+ \frac{ (a + 1)\lambda^2}{2}I\bigl(\llvert\beta\rrvert> a \lambda
\bigr)\qquad\mbox{for some } a > 2,
\end{eqnarray*}
and for the MCP penalty function,
\[
p_\lambda\bigl(\llvert\beta\rrvert\bigr) = \lambda\biggl( \llvert
\beta\rrvert- \frac{\beta
^2}{2a\lambda} \biggr) I\bigl(0 \leq\llvert\beta\rrvert< a
\lambda\bigr) + \frac
{a\lambda^2}{2}I \bigl(\llvert\beta\rrvert\geq a \lambda
\bigr)\qquad\mbox{for some } a > 1.
\]
For both penalty functions, the tuning parameter $\lambda$ controls
the complexity
of the selected model
and goes to zero as $n$ increases to $\infty$.

\subsection{Solving the penalized estimator}\label{sec3.2}
We propose an effective algorithm to solve the above penalized
estimation problem. The algorithm is largely based on
the idea of the local linear approximation (LLA) [\citet{lla}].
We employ a new trick based on the observation $\llvert \beta
_j\rrvert =\rho_\tau
(\beta_j) + \rho_\tau(-\beta_j)$
to transform the approximated objective function to a quantile
regression objective function based
on an augmented data set, so that the penalized estimator can be
obtained by iteratively
solving unpenalized weighted quantile regression problems.

More specifically, we initialize the algorithm by setting $\vbeta=0$
and $\vxi=0$. Then
for each step $t \geq1$, we update the estimator by
%
%
\begin{equation}\label{q_scad_step}
\bigl(\hat{\vbeta}^t, \hat{\vxi}^t \bigr)
= \mathop{\operatorname{argmin}}_{(\vbeta, \vxi)} \Biggl\{n^{-1}\sum
_{i=1}^n\rho_\tau
\bigl(Y_i - \vx_i'\vbeta- \bolds{\Pi}(
\vz_i)'\vxi\bigr)+ \sum_{j=1}^{p_n}
p'_\lambda\bigl(\bigl\llvert\hat{\beta}_j^{t-1}
\bigr\rrvert\bigr)\llvert\beta_j\rrvert\Biggr\},\hspace*{-30pt}
\end{equation}
where $\hat{\beta}_j^{t-1}$ is the value of $\beta_j$ at step $t-1$.\vspace*{1pt}

By observing that we can write $\llvert \beta_j\rrvert $ as $\rho_\tau
(\beta_j) +
\rho_\tau(-\beta_j)$,
the above minimization problem
can be framed as an unpenalized weighted quantile regression problem
with $n+2p_n$
augmented observations. We denote these
augmented observations by $(Y_i^*, \vx_i^*, \vz_i^*)$, $i=1,\ldots,
(n+2p_n)$.
The first $n$ observations are those in the original data, that is,
$(Y_i^*, \vx_i^*, \vz_i^*)=(Y_i, \vx_i, \vz_i)$,
$i=1,\ldots,n$;
for the next $p_n$ observations, we have $(Y_i^*, \vx_i^*, \vz_i
^*)=(0,1,0)$, $i=n+1,\ldots, n+p_n$;
and the last $p_n$ observations are given by $(Y_i^*, \vx_i^*, \vz_i
^*)=(0,-1,0)$,
$i=n+p_n+1,\ldots, n+2p_n$. We fit weighted
linear quantile regression model with the observations $(Y_i^*, \vx_i
^*, \vz_i^*)$
and corresponding weights $w_i^{t*}$, where $w_i^{t*}=1$, $i=1,\ldots,n$;
$w_{n+j}^{t*}=p'_\lambda(\llvert \hat{\beta}_j^{t-1}\rrvert )$,
$j=1,\ldots,p_n$;
and
$w_{n+p_n+j}^{t*}=-p'_\lambda(\llvert \hat{\beta}_j^{t-1}\rrvert )$,
$j=1,\ldots,p_n$.\vspace*{2pt}

The above new algorithm is simple and convenient, as
weighted quantile regression can be implemented using many existing
software packages.
In our simulations, we used the quantreg package in R and
continue with the iterative procedure until $\llVert \hat{\vbeta}^t -
\hat
{\vbeta}^{t-1}\rrVert _1 < 10^{-7}$.

\subsection{Asymptotic theory}\label{sec3.3}

In addition to Conditions \ref{cond_f}--\ref{cond_sigma_large_p},
we impose an additional condition on how quickly a nonzero signal can
decay, which is
needed to identify the underlying model.

%
\begin{condition}[(Condition on the signal)]\label{cond_small_sig}
There exist positive constants $C_4$ and
$C_5$ such that $2C_3 < C_4 < 1$ and $n^{(1-C_4)/2} \mathop{\min}_{1
\leq j \leq q_n} \llvert \beta_{0j}\rrvert \geq C_5$.
\end{condition}

Due to the nonsmoothness and nonconvexity of the penalized objective
function $Q^P(\vbeta,\vxi)$,
the classical KKT condition is not applicable to analyzing the
asymptotic properties of the penalized estimator.
To investigate the asymptotic theory of the nonconvex estimator
for ultra-high dimensional partially linear additive quantile
regression model,
we explore the necessary condition for the local minimizer
of a convex differencing problem [\citet{TaoAn97};
\citet{quantUltraHigh}] and extend it to the setting involving
nonparametric
components.


Our approach concerns a nonconvex
objective function that can be expressed as the difference of two
convex functions. Specifically, we consider objective functions
belonging to the class
\[
\mathbf{F} = \bigl\{ q(\veta): q(\veta) = k(\veta) - l(\veta), k(\cdot),l(
\cdot) \mbox{ are both convex} \bigr\}.
\]
This is a very general formulation that incorporates many different forms
of penalized objective functions.
The subdifferential of $k(\veta)$ at $\veta=\veta_0$ is defined as
\[
\partial k(\veta_0) = \bigl\{t: k(\veta) \geq k(
\veta_0) + (\veta-\veta_0)'t, \forall
\veta\bigr\}.
\]
Similarly, we can define the subdifferential of $l(\veta)$.
Let $\operatorname{dom}(k) =
\{ \veta: k(\veta) < \infty\}$ be the effective domain
of $k$.
A necessary condition for $\veta^*$ to be a local minimizer of
$q(\veta)$ is that
$\veta^*$ has a neighborhood $U$ such that
$\partial l(\veta) \cap\partial k(\veta^*) \neq\varnothing, \forall
\veta\in U \cap\operatorname{dom}(k)$
(see Lemma \ref{lem_diff_convex} in the \hyperref[appe]{Appendix}).


To appeal to the above necessary condition for the convex differencing problem,
it is noted that $Q^P(\vbeta, \vxi)$ can be written as
\[
Q^P(\vbeta,\vxi) = k(\vbeta,\vxi) - l(\vbeta, \vxi),
\]
where the two convex functions $k(\vbeta,\vxi) =n^{-1} \sum_{i=1}^n\rho
_\tau(Y_i - \vx_i'\vbeta-\bolds{\Pi}(\vz_i)'\vxi) + \lambda\sum_{j=1}^{p_n}
\llvert \beta_j\rrvert $, and $l(\vbeta, \vxi) = \sum_{j=1}^{p_n}
L(\beta_j)$.
The specific form of
$L(\beta_j)$ depends on the penalty function being used.
For the SCAD penalty function,
\begin{eqnarray*}
L(\beta_j) &=& \bigl[ \bigl(\beta_j^2 + 2
\lambda\llvert\beta_j\rrvert+ \lambda^2 \bigr)/
\bigl(2(a-1)\bigr) \bigr]I\bigl(\lambda\leq\llvert\beta_j\rrvert\leq
a\lambda\bigr)
\\
&&{}+ \bigl[ \lambda\llvert\beta_j\rrvert- (a+1)
\lambda^2/2 \bigr]I\bigl(\llvert\beta_j\rrvert> a\lambda
\bigr);
\end{eqnarray*}
while for the MCP penalty function,
\begin{eqnarray*}
&& L(\beta_j) = \bigl[ \beta_j^2 / (2a) \bigr]I
\bigl(0 \leq\llvert\beta_j\rrvert< a\lambda\bigr) + \bigl[
\lambda\llvert\beta_j\rrvert- a \lambda^2 / 2 \bigr]I
\bigl(\llvert\beta_j\rrvert\geq a \lambda\bigr).
\end{eqnarray*}


Building on the convex differencing structure, we show that
with probability approaching one that the oracle estimator
$(\hat{\vbeta}', \hat{\vxi}')'$, where
$\hat{\vbeta}= (\hat{\vbeta}_1', \mathbf{0}_{p_n-q_n}' )'$,
is a local minimizer of $Q^P(\vbeta,\vxi)$.
To study the necessary optimality condition, we formally define
$\partial k(\vbeta,\vxi)$ and $\partial l(\vbeta, \vxi)$, the
subdifferentials of $k(\vbeta, \vxi)$ and $l(\vbeta, \vxi)$, respectively.
First, the function $l(\vbeta, \vxi)$ does not depend on $\vxi$ and
is differentiable everywhere.
Hence, its subdifferential is simply the regular derivative. For any
value of $\vbeta$ and $\vxi$,
\begin{eqnarray*}
\partial l(\vbeta, \vxi) &=& \biggl\{ \mu = (\mu_1,
\mu_2,\ldots,\mu_{p_n+L_n})' \in
\mathds{R}^{p_n+L_n}:
\\
&&{} \mu_j = \frac{\partial l(\vbeta)}{\partial\beta_j}, 1 \leq j \leq p_n;
\mu_j=0, p_n+1 \leq j \leq p_n+L_n
\biggr\}.
\end{eqnarray*}
%
For $1 \leq j \leq p_n$, for the SCAD penalty function,
\[
\frac{\partial l(\vbeta)}{\partial\beta_j} = \cases{ 0, &\quad$0 \leq
\llvert\beta_j\rrvert<
\lambda$,
\vspace*{3pt}\cr
\bigl(\beta_j-\lambda\operatorname{sgn}(
\beta_j)\bigr)/(a-1), &\quad$\lambda\leq\llvert\beta_j
\rrvert\leq a\lambda$,
\vspace*{3pt}\cr
\lambda\operatorname{sgn}(\beta_j), &\quad$
\llvert\beta_j\rrvert>a\lambda$,}
\]
while for the MCP penalty function,
\[
\frac{\partial l(\vbeta)}{\partial\beta_j} = \cases{ \beta_j/a, &\quad
$0 \leq\llvert
\beta_j\rrvert< a\lambda$,
\vspace*{3pt}\cr
\lambda\operatorname{sgn}(
\beta_j), &\quad$\llvert\beta_j\rrvert\geq a\lambda$.}
\]
On the other hand, the function $k(\vbeta, \vxi)$ is not
differentiable everywhere.
Its subdifferential at $(\vbeta,\vxi)$ is a collection of $(p_n+L_n)$-vectors:
\begin{eqnarray*}
\partial k(\vbeta, \vxi) &=& \Biggl\{ \vkappa= ( \kappa_1, \kappa
_2,\ldots,\kappa_{p_n+L_n})' \in
\mathds{R}^{p_n+L_n}:
\\
&&{} \kappa_j = -\tau n^{-1}\sum
_{i=1}^nx_{ij}I\bigl(Y_i -
\vx_i'\vbeta- \bolds{\Pi}(\vz_i)
'\vxi> 0\bigr)
\\
&&{} + (1-\tau) n^{-1}\sum_{i=1}^nx_{ij}I
\bigl(Y_i - \vx_i'\vbeta- \bolds{\Pi}(
\vz_i)'\vxi< 0\bigr)
\\
&&{}- n^{-1}\sum_{i=1}^nx_{ij}a_i
+ \lambda l_j,\mbox{ for } 1 \leq j \leq p_n;
\\
&&{} \kappa_j = -\tau n^{-1}\sum
_{i=1}^n\Pi_{j-p_n}(\vz_i)I
\bigl(Y_i - \vx_i '\vbeta- \bolds{\Pi}(
\vz_i)'\vxi> 0\bigr)
\\
&&{} +(1-\tau)n^{-1}\sum_{i=1}^n
\Pi_{j-p_n}(\vz_i) I\bigl(Y_i - \vx
_i'\vbeta- \bolds{\Pi}(\vz_i)'
\vxi< 0\bigr)
\\
&&{} -n^{-1}\sum_{i=1}^n
\Pi_{j-p_n}(\vz_i)a_i, \mbox{ for }
p_n + 1 \leq j \leq p_n +L_n \Biggr\},
\end{eqnarray*}
where we write $\bolds{\Pi}(\vz_i)=(1, \Pi_1(\vz_i), \ldots, \Pi
_{L_n}(\vz_i))'$;
$a_i=0$ if $Y_i-\vx_i'\vbeta- \bolds{\Pi}(\vz_i)'\vxi\neq0$ and
$a_i \in[\tau
-1,\tau]$ otherwise;
for $1 \leq j \leq p_n$,
$l_j = \operatorname{sgn}(\beta_j)$ if $\beta_j \neq0$ and $l_j \in[-1,1]$
otherwise.

In the following, we analyze the subgradient of the unpenalized
objective function, which plays an essential
role in checking the condition of the optimality condition. The
subgradient $s (\vbeta,\vxi) = (s_1(\vbeta,\vxi
),\ldots,s_{p_n}(\vbeta,\vxi),\ldots,\break s_{p_n+L_n}(\vbeta,\vxi
) )'$ is given by
\begin{eqnarray*}
s_j(\vbeta, \vxi) &=& -\frac{\tau}{n}\sum
_{i=1}^nx_{ij}I\bigl(Y_i -
\vx_i'\vbeta- \bolds{\Pi}(\vz_i)'
\vxi> 0\bigr)
\\
&&{} + \frac{1-\tau}{n} \sum_{i=1}^nx_{ij}I
\bigl(Y_i - \vx_i'\vbeta- \bolds{\Pi}(
\vz_i) '\vxi< 0\bigr)
\\
&&{}- \frac{1}{n} \sum_{i=1}^nx_{ij}a_i\qquad\mbox{for } 1 \leq j \leq p_n,
\\
s_j(\vbeta, \vxi) &=& -\frac{\tau}{n}\sum
_{i=1}^n\Pi_{j-p_n}(\vz_i)I
\bigl(Y_i - \vx_i'\vbeta- \bolds{\Pi}(
\vz_i)'\vxi> 0\bigr)
\\
&&{}+ \frac{1-\tau}{n} \sum_{i=1}^n
\Pi_{j-p_n}(\vz_i)I\bigl(Y_i - \vx
_i'\vbeta- \bolds{\Pi}(\vz_i)'
\vxi< 0\bigr)
\\
&&{}- \frac{1}{n} \sum_{i=1}^n\Pi
_{j-p_n}(\vz_i)a_i\qquad\mbox{for } p_n
+ 1 \leq j \leq p_n+L_n,
\end{eqnarray*}
where $a_i$ is defined as before.
The following lemma states the behavior of $s_j(\hat{\vbeta},\hat
{\vxi})$ when being evaluated
at the oracle estimator.

%
\begin{lemma}
\label{lem_sub_diff_q}
Assume Conditions \ref{cond_f}--\ref{cond_small_sig} are satisfied,
$\lambda= o (n^{-(1-C_4)/2} )$,\break  $n^{-1/2}q_n = o(\lambda
)$, $n^{-1/2}k_n = o(\lambda)$ and $\log(p_n) = o(n\lambda^2)$. For
the oracle estimator $ (\hat{\vbeta},\hat{\vxi} )$ there
exists $a_i^*$ with $a_i^*=0$ if $Y_i - \vx_i'\hat{\vbeta} - \bolds
{\Pi}(\vz_i)
'\hat{\vxi}\neq0$ and $a_i^* \in[\tau-1,\tau]$ otherwise, such
that for $s_j(\hat{\vbeta},\hat{\vxi})$ with $a_i=a_i^*$,
with probability approaching one
%
%
\begin{eqnarray}
s_j (\hat{\vbeta},\hat{\vxi} ) &=&0,\qquad j=1,\ldots,q_n
\mbox{ or } j = p_n+1,\ldots,p_n+L_n,
\label
{lem_sub_diff_q_1}
\\
\llvert\hat{\beta}_j\rrvert&\geq& (a+1/2)\lambda,\qquad j=1,\ldots
,q_n, \label{lem_sub_diff_q_2}
\\
\bigl\llvert s_j (\hat{\vbeta},\hat{\vxi} )\bigr\rrvert&\leq& c
\lambda\qquad\forall c>0, j=q_n+1,\ldots,p_n. \label{lem_sub_diff_q_3}
\end{eqnarray}
\end{lemma}

\begin{rem*}
Note that for $\kappa_{j} \in\partial k(\vbeta,\vxi)$ and $l_j$ as
defined earlier
\begin{eqnarray*}
\kappa_j &=& s_j(\vbeta,\vxi) + \lambda
l_j\qquad\mbox{for } 1 \leq j \leq p_n
\quad\mbox{and}
\\
\kappa_{j} &=& s_{j}(\vbeta,\vxi)\qquad\mbox{for }
p_n+1 \leq j \leq p_n+L_n.
\end{eqnarray*}


Thus, Lemma \ref{lem_sub_diff_q} provides important insight on the
asymptotic behavior of $\vkappa\in\partial k(\hat{\vbeta},\hat
{\vxi})$. Consider a small neighborhood around the oracle estimator
$ (\hat{\vbeta},\hat{\vxi} )$ with radius $\lambda/ 2$.
Building on Lemma \ref{lem_sub_diff_q}, we prove in the \hyperref[appe]{Appendix} that
with probability tending to one, for any $ (\vbeta, \vxi)
\in\mathds{R}^{p_n+L_n}$
in this neighborhood, there exists $\vkappa= (\kappa_1,\ldots
,\kappa_{p_n}, \mathbf{0}_{L_n}' )' \in\partial k(\hat
{\vbeta},\hat{\vxi})$ such that
\begin{eqnarray*}
\frac{\partial l(\vbeta, \vxi)}{\partial\beta_j} &=& \kappa_j,\qquad j=1,\ldots
,p_n\quad\mbox{and}
\\
\frac{\partial l(\vbeta, \vxi)}{\partial\vxi_j} &=& \kappa_{p_n+j},\qquad
j=1,\ldots,L_n.
\end{eqnarray*}
This leads to the main theorem of the paper.
Let $\mathcal{E}_n(\lambda)$ be the set of local minima of
$Q^P(\vbeta,\vxi)$. The theorem below shows that with probability
approaching one, the oracle estimator belongs to the set $\mathcal
{E}_n(\lambda)$.
\end{rem*}

%
\begin{theorem}
\label{scad_local_min}
Assume Conditions \ref{cond_f}--\ref{cond_small_sig} are satisfied.
Consider either the SCAD or the MCP penalty function with tuning
parameter $\lambda$. Let $\hat{\veta} \equiv(\hat{\vbeta},
\hat{\vxi} )$ be the oracle estimator. If 
$\lambda= o (n^{-(1-C_4)/2} ), n^{-1/2}q_n = o(\lambda)$,
$n^{-1/2}k_n = o(\lambda)$ and $\log(p_n) = o(n\lambda^2)$,
then
\[
P \bigl(\hat{\veta} \in\mathcal{E}_n(\lambda) \bigr) \rightarrow1
\qquad\mbox{as } n\rightarrow\infty.
\]
\end{theorem}

\begin{rem*}
The conditions for $\lambda$ in the theorem are satisfied for $\lambda
= n^{-1/2+\delta}$ where $\delta\in( \max(1/(2r+1),C_3), C_4 )$. The
fastest rate of $p_n$ allowed is
$p_n = \operatorname{exp}(n^{\alpha})$ with $0 < \alpha< 1/2 + 2\delta$. Hence,
we allow for the ultra-high dimensional setting.
\end{rem*}

\begin{rem*}
The selection of the tuning parameter
$\lambda$ is important in practice.
Cross-validation is a common approach, but is known to often result in
overfitting.
\citet{Lee}
recently proposed high dimensional BIC for linear quantile
regression when $p$ is much larger than $n$. Motivated
by their work, we
choose $\lambda$ that minimizes the following high dimensional BIC criterion:
%
%
\begin{eqnarray}\label{large_p_bic}
\mbox{QBIC}(\lambda) &=& \log\Biggl( \sum
_{i=1}^n\rho_\tau\bigl(Y_i
- \vx_i'\hat{\vbeta}_{\lambda} - \bolds{\Pi}(
\vz_i)'\hat{\vxi}_{\lambda
} \bigr) \Biggr)
\nonumber\\[-8pt]\\[-8pt]\nonumber
&&{} +
\nu_{\lambda} \frac{ \log(p_n)\log(\log
(n))}{2n}, 
\end{eqnarray}
where $p_n$ is the number of candidate linear covariates and $\nu
_\lambda$ is the degrees of freedom of the
fitted model, which is the number of interpolated fits for quantile regression.
\end{rem*}

\section{Simulation}\label{sec4}
We investigate the performance of the penalized partially linear additive
quantile regression estimator
in high dimension. We focus on the SCAD penalty and referred to the new
procedure as Q-SCAD.
An alternative popular nonconvex penalty function is
the MCP penalty [\citet{Zhang}], the simulation results for which are
found to be similar and reported in the online supplementary material
[\citet{Supp}].
The Q-SCAD is compared with three alternative procedures:
partially linear additive quantile regression estimator with the LASSO
penalty (Q-LASSO),
partially linear additive mean regression with SCAD penalty (LS-SCAD) and
LASSO penalty (LS-LASSO). It worth noting that for the mean
regression case,
there appears to be no theory in the literature for the ultra-high
dimensional case.

We first generate $\tilde{X}=(\tilde{X}_1,\ldots, \tilde{X}_{p+2})'$
from the $N_{p+2}(\mathbf{0}_{p+2},\Sigma)$ multivariate normal distribution,
where $\Sigma= (\sigma_{jk})_{ (p+2) \times(p+2)}$ with $\sigma
_{jk} = 0.5^{|j-k|}$. Then we set $X_1 = \sqrt{12}\Phi(\tilde{X}_1)$
where $\Phi(\cdot)$ is
distribution function of $N(0,1)$ distribution
and $\sqrt{12}$ scales $X_1$ to have standard deviation one.
Furthermore, we let $Z_{1} = \Phi(\tilde{X}_{25})$, $Z_2 = \Phi
(\tilde{X}_{26})$, $X_i = \tilde{X}_i$ for $i=2,\ldots,24$ and $X_i
= \tilde{X}_{i-2}$ for $i = 27,\ldots,p+2$.
The random responses are generated from the regression model
%
%
\begin{equation}
Y_i = X_{i6}\beta_1 + X_{i12}
\beta_2 + X_{i15}\beta_3 + X_{i20}
\beta_4 + \sin(2 \pi Z_{i1}) + Z_{i2}^3
+ \varepsilon_i, 
\end{equation}
where $\beta_j \sim U[0.5,1.5]$ for $1 \leq j \leq4$. We consider
three different distributions of the error term $\varepsilon_i$: (1)
standard normal distribution; (2) $t$ distribution with 3 degrees of
freedom; and (3) heteroscedastic
normal distribution $\varepsilon_i = \tilde{X}_{i1}\zeta_i$ where
$\zeta_i \sim N(0,\sigma=0.7)$
are independent of the $X_i$'s. 

We perform 100 simulations for each setting with sample size $n=300$,
and $p=100$, $300$, $600$. 
Results for additional simulations with sample sizes of $50$, 100 and
200 are provided in the online supplementary material [\citet
{Supp}]. For the heteroscedastic error case, we model $\tau= 0.7$
and $0.9$; otherwise, we model the conditional median.
Note that at $\tau= 0.7$ or 0.9, when the error has the
aforementioned heteroscedastic distribution, $X_1$ is part of the true
model. At these two quantiles, the true model consists of 5 linear
covariates. In all simulations, the number of basis functions is set to three,
which we find to work satisfactorily in a variety of settings. For the
LASSO method, we select the tuning parameters $\lambda$
by using five-fold cross validation. For the Q-SCAD model, we select
$\lambda$ that minimizes (\ref{large_p_bic}) while for LS-SCAD we use
a least squares equivalent. The tuning parameter $a$ in the SCAD
penalty function
is set to 3.7 as recommended in \citet{fanLi}. To assess the
performance of different methods, we adopt the following criteria:
\begin{longlist}
%
\item[1.] False Variables (FV): average number of nonzero linear covariates
incorrectly included in the model.
\item[2.] True Variables (TV): 
average number of nonzero linear covariates correctly included in the model.
\item[3.] True: proportion of times the true model is exactly identified.
\item[4.] P: proportion of times $X_{1}$ is selected.
\item[5.] AADE: average of the \textit{average absolute deviation} (ADE) of
the fit of the nonlinear components,
where the ADE is defined as $n^{-1}\sum_{i=1}^n|\hat{g}(\vz_i) -
g_0(\vz_i)|$.
\item[6.] MSE: average of the mean squared error for estimating $\vbeta
_0$, that is, the average of $\llVert \hat{\vbeta}-\vbeta_0\rrVert ^2$
across all simulation runs.
\end{longlist}




%
%
\begin{table}[t]
\tabcolsep=0pt
\caption{Simulation results comparing quantile ($\tau=0.5$) and mean regression
using SCAD and LASSO penalty functions for $\varepsilon\sim N(0,1)$}\label{tab1}
\begin{tabular*}{\tablewidth}{@{\extracolsep{\fill}}@{}lccd{2.2}ccccc@{}}
\hline
\textbf{Method} & $\bolds{n}$ & $\bolds{p}$ & \multicolumn{1}{c}{\textbf{FV}} & \textbf{TV} & \textbf{True} & \textbf{P} & \textbf{AADE} & \textbf{MSE}\\
\hline
Q-SCAD & 300 & 100 & 0.20 & 4.00 & 0.88 & 0.00 & 0.16 & 0.03 \\
Q-LASSO & 300 & 100 & 12.88 & 4.00 & 0.00 & 0.13 & 0.16 & 0.13 \\
LS-SCAD & 300 & 100 & 0.32 & 4.00 & 0.85 & 0.00 & 0.13 & 0.02 \\
LS-LASSO & 300 & 100 & 11.63 & 4.00 & 0.00 & 0.12 & 0.13 & 0.07
\\[3pt]
Q-SCAD & 300 & 300 & 0.04 & 4.00 & 0.96 & 0.00 & 0.15 & 0.02 \\
Q-LASSO & 300 & 300 & 15.93 & 4.00 & 0.00 & 0.07 & 0.16 & 0.14 \\
LS-SCAD & 300 & 300 & 0.33 & 4.00 & 0.78 & 0.00 & 0.12 & 0.02 \\
LS-LASSO & 300 & 300 & 15.00 & 4.00 & 0.00 & 0.04 & 0.13 & 0.09
\\[3pt]
Q-SCAD & 300 & 600 & 0.06 & 4.00 & 0.94 & 0.00 & 0.15 & 0.02 \\
Q-LASSO & 300 & 600 & 21.86 & 4.00 & 0.01 & 0.06 & 0.16 & 0.16 \\
LS-SCAD & 300 & 600 & 2.57 & 4.00 & 0.69 & 0.01 & 0.13 & 0.06 \\
LS-LASSO & 300 & 600 & 17.11 & 4.00 & 0.00 & 0.04 & 0.13 & 0.09 \\
\hline
\end{tabular*}
\end{table}

%
%
\begin{table}[b]
\tabcolsep=0pt
\caption{Simulation results comparing quantile ($\tau=0.5$) and mean regression
using SCAD and LASSO penalty functions for $\varepsilon\sim T_3$}\label{tab2}
\begin{tabular*}{\tablewidth}{@{\extracolsep{\fill}}@{}lccd{2.2}ccccc@{}}
\hline
\textbf{Method} & $\bolds{n}$ & $\bolds{p}$ & \multicolumn{1}{c}{\textbf{FV}} & \textbf{TV} & \textbf{True} & \textbf{P} & \textbf{AADE} & \textbf{MSE}\\
\hline
Q-SCAD & 300 & 100 & 0.07 & 4.00 & 0.95 & 0.00 & 0.16 & 0.03 \\
Q-LASSO & 300 & 100 & 13.09 & 4.00 & 0.01 & 0.17 & 0.17 & 0.15 \\
LS-SCAD & 300 & 100 & 1.08 & 3.99 & 0.45 & 0.02 & 0.19 & 0.11 \\
LS-LASSO & 300 & 100 & 10.15 & 3.94 & 0.02 & 0.08 & 0.19 & 0.31
\\[3pt]
Q-SCAD & 300 & 300 & 0.05 & 4.00 & 0.97 & 0.00 & 0.17 & 0.03 \\
Q-LASSO & 300 & 300 & 18.42 & 4.00 & 0.00 & 0.08 & 0.18 & 0.18 \\
LS-SCAD & 300 & 300 & 1.22 & 4.00 & 0.46 & 0.00 & 0.20 & 0.11 \\
LS-LASSO & 300 & 300 & 15.15 & 3.99 & 0.01 & 0.08 & 0.21 & 0.26
\\[3pt]
Q-SCAD & 300 & 600 & 0.06 & 3.98 & 0.94 & 0.00 & 0.16 & 0.04 \\
Q-LASSO & 300 & 600 & 20.81 & 4.00 & 0.01 & 0.03 & 0.18 & 0.23 \\
LS-SCAD & 300 & 600 & 1.33 & 4.00 & 0.45 & 0.00 & 0.19 & 0.14 \\
LS-LASSO & 300 & 600 & 17.40 & 4.00 & 0.01 & 0.01 & 0.20 & 0.28 \\
\hline
\end{tabular*}
\end{table}

%
%
\begin{table}[t]
\tabcolsep=0pt
\caption{Simulation results comparing quantile ($\tau=0.7$) and mean regression
using SCAD and LASSO penalty functions for heteroscedastic errors}\label{tab3}
\begin{tabular*}{\tablewidth}{@{\extracolsep{\fill}}@{}lccd{2.2}ccccc@{}}
\hline
\textbf{Method} & $\bolds{n}$ & $\bolds{p}$ & \multicolumn{1}{c}{\textbf{FV}} & \textbf{TV} & \textbf{True} & \textbf{P} & \textbf{AADE} & \textbf{MSE}\\
\hline
Q-SCAD & 300 & 100 & 0.21 & 4.84 & 0.70 & 0.84 & 0.17 & 0.05 \\
Q-LASSO & 300 & 100 & 13.86 & 4.97 & 0.00 & 0.97 & 0.24 & 0.15 \\
LS-SCAD & 300 & 100 & 1.09 & 4.06 & 0.01 & 0.06 & 0.16 & 0.69 \\
LS-LASSO & 300 & 100 & 11.48 & 4.13 & 0.00 & 0.13 & 0.17 & 0.78
\\[3pt]
Q-SCAD & 300 & 300 & 0.20 & 4.77 & 0.61 & 0.77 & 0.20 & 0.06 \\
Q-LASSO & 300 & 300 & 18.54 & 4.97 & 0.00 & 0.97 & 0.27 & 0.18 \\
LS-SCAD & 300 & 300 & 3.28 & 4.00 & 0.00 & 0.00 & 0.16 & 0.68 \\
LS-LASSO & 300 & 300 & 15.85 & 4.08 & 0.00 & 0.08 & 0.16 & 0.79
\\[3pt]
Q-SCAD & 300 & 600 & 0.16 & 4.59 & 0.48 & 0.59 & 0.26 & 0.08 \\
Q-LASSO & 300 & 600 & 23.26 & 4.89 & 0.00 & 0.89 & 0.31 & 0.24 \\
L
LS-SCAD & 300 & 600 & 6.31 & 4.02 & 0.00 & 0.02 & 0.16 & 0.69 \\
LS-LASSO & 300 & 600 & 18.50 & 4.09 & 0.00 & 0.09 & 0.16 & 0.83 \\
\hline
\end{tabular*}
\end{table}

%
%
\begin{table}[b]
\tabcolsep=0pt
\caption{Simulation results comparing quantile ($\tau=0.9$) and mean regression
using SCAD and LASSO penalty functions for heteroscedastic errors}\label{tab4}
\begin{tabular*}{\tablewidth}{@{\extracolsep{\fill}}@{}lccd{2.2}ccccc@{}}
\hline
\textbf{Method} & $\bolds{n}$ & $\bolds{p}$ & \multicolumn{1}{c}{\textbf{FV}} & \textbf{TV} & \textbf{True} & \textbf{P} & \textbf{AADE} & \textbf{MSE}\\
\hline
Q-SCAD & 300 & 100 & 0.06 & 4.93 & 0.91 & 0.98 & 0.24 & 0.30 \\
Q-LASSO & 300 & 100 & 12.94 & 5.00 & 0.00 & 1.00 & 0.49 & 0.73 \\
LS-SCAD & 300 & 100 & 1.09 & 4.06 & 0.01 & 0.06 & 0.16 & 4.72 \\
LS-LASSO & 300 & 100 & 11.48 & 4.13 & 0.00 & 0.13 & 0.17 & 4.73
\\[3pt]
Q-SCAD & 300 & 300 & 0.26 & 5.00 & 0.81 & 1.00 & 0.19 & 0.24 \\
Q-LASSO & 300 & 300 & 16.33 & 5.00 & 0.00 & 1.00 & 0.62 & 0.92 \\
LS-SCAD & 300 & 300 & 3.28 & 4.00 & 0.00 & 0.00 & 0.16 & 4.63 \\
LS-LASSO & 300 & 300 & 15.85 & 4.08 & 0.00 & 0.08 & 0.16 & 4.67
\\[3pt]
Q-SCAD & 300 & 600 & 0.34 & 4.94 & 0.77 & 1.00 & 0.21 & 0.29 \\
Q-LASSO & 300 & 600 & 19.79 & 4.97 & 0.00 & 1.00 & 0.74 & 1.15 \\
LS-SCAD & 300 & 600 & 6.31 & 4.02 & 0.00 & 0.02 & 0.16 & 4.64 \\
LS-LASSO & 300 & 600 & 18.50 & 4.09 & 0.00 & 0.09 & 0.16 & 4.74 \\
\hline
\end{tabular*}
\end{table}

The simulation results are summarized in Tables~\ref{tab1}--\ref{tab4}.
Tables \ref{tab1}~and~\ref{tab2} correspond to $\tau= 0.5$, $N(0,1)$ and $T_3$ error
distribution, respectively.
Tables \ref{tab3}~and~\ref{tab4} are for the heteroscedastic error, $\tau= 0.7$ and
$0.9$, respectively. Least squares based estimates of $\hat{\vbeta}$
for $\tau=0.7$ or $0.9$ are obtained by assuming $\varepsilon_i \sim
N(0,\sigma)$, with estimates of $\sigma$ being used in each
simulation. An extension of Table~\ref{tab3} for $p = 1200$ and 2400 is included
in the online supplementary material [\citet{Supp}].
We observe that the method with the SCAD penalty tends to pick a
smaller and more accurate model. The advantages of quantile regression
can be seen by its stronger performance
at the presence of heavy-tailed distribution or heteroscedastic errors.
For the latter case,
the least squared based methods perform poorly in identifying the
active variables in the
dispersion function. 
Estimation of the nonlinear terms is similar across different error
distributions and
different values of~$p$.

\section{An application to birth weight data}\label{sec5}
Votavova et~al. (\citeyear{AppliedExample}) collected blood samples
from peripheral blood,
cord blood and the placenta from 20 pregnant smokers and 52 pregnant
women without significant exposure to smoking. Their main objective was
to identify the difference in transcriptome alterations between the two
groups. Birth weight of the baby (in kilograms) was recorded along with
age of the mother, gestational age, parity, measurement of the amount
of cotinine, a chemical found in tobacco, in the blood and mother's
BMI. Low birth weight is known to be associated with both short-term
and long-term health complications. Scientists are interested in which
genes are associated with low birth weight
[\citet{birthWt2}].

We consider modeling the 0.1, 0.3 and 0.5 conditional quantiles of
infant birth weight.
We use the genetic data from the peripheral blood sample which include
64 subjects after dropping those with incomplete information. The blood
samples were assayed using HumanRef-8 v3 Expression BeadChips with
24,539 probes. For each quantile, the top 200 probes are selected using
the quantile-adaptive screening method [\citet{heWang}]. The
gene expression values of the 200 probes are included as linear
covariates for the semiparametric quantile regression model. The
clinical variables parity, gestational age, cotinine level and BMI
are also included as linear covariates. The age of the mother is
modeled nonparametrically as exploratory analysis reveals potential
nonlinear effect.

We consider the semiparametric quantile regression model with the SCAD
and LASSO penalty functions.
Least squares based semiparametric models with the SCAD and LASSO
penalty functions are also considered. Results for the MCP penalty are
reported in the online supplementary material [\citet{Supp}].
The tuning parameter $\lambda$ is selected by minimizing (\ref
{large_p_bic}) for the SCAD estimator and by five-fold cross validation
for LASSO as discussed in Section~\ref{sec4}. The third column of
Table~\ref{applied_models} reports the number of nonzero elements, ``Original
NZ,'' for each model. As expected, the LASSO method selects a larger
model than the SCAD penalty does. The number of nonzero variables
varies with the quantile level, providing evidence that mean regression
alone would provide a limited view of the conditional distribution.

%
\begin{table}[b]
\tabcolsep=0pt

\caption{Quantile ($\tau=0.1$, 0.3 and 0.5) and mean regression
analysis of birth weight based on the original data and the random partitioned data}\label{tab5}\label{applied_models}
\begin{tabular*}{\tablewidth}{@{\extracolsep{\fill}}@{}lcd{2.0}cd{2.2}@{}}
\hline
$\bolds{\tau}$ & \textbf{Method} & \multicolumn{1}{c}{\textbf{Original NZ}} & \textbf{Prediction error} & \multicolumn{1}{c@{}}{\textbf{Randomized NZ}}\\
\hline
0.10 & Q-SCAD & 2 & 0.07 (0.03) & 2.27 \\
0.10 & Q-LASSO & 10 & 0.08 (0.02) & 3.09
\\[3pt]
0.30 & Q-SCAD & 7 & 0.18 (0.04) & 6.74 \\
0.30 & Q-LASSO & 22 & 0.16 (0.03) & 12.39
\\[3pt]
0.50 & Q-SCAD & 5 & 0.21 (0.04) & 5.80 \\
0.50 & Q-LASSO & 6 & 0.20 (0.04) & 14.25
\\[3pt]
Mean & LS-SCAD &12 & 0.20 (0.04) & 5.43 \\
Mean & LS-LASSO & 12 & 0.20 (0.04) & 3.77 \\
\hline
\end{tabular*}
\end{table}

Next, we compare different models on 100 random partitions of the data set.
For each partition, we randomly select 50 subjects for the training
data and 14 subjects for the test data.
The fourth column of Table~\ref{applied_models} reports the prediction
error evaluated on the test data, defined as
$14^{-1}\sum_{i=1}^{14} \rho_\tau(Y_i - \hat{Y}_i)$; while the
fifth column reports the average number of linear covariates included
in each model (denoted by ``Randomized NZ''). Standard errors for the
prediction error is reported in parentheses. We note that the SCAD
method produces notably smaller models than the Lasso method does
without sacrificing much prediction accuracy.

%
\begin{figure}[b]

\includegraphics{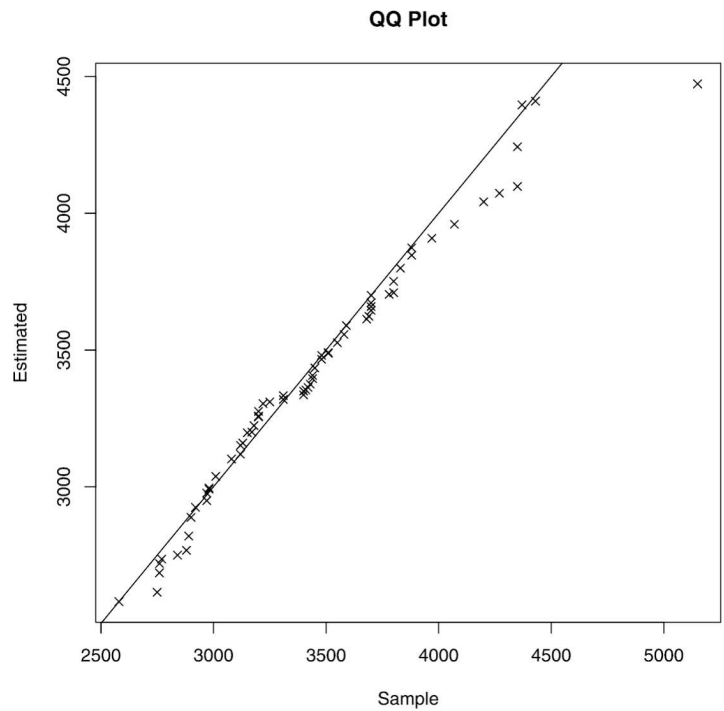}

\caption{Lack-of-fit diagnostic QQ plot for the birth weight data
example.}\label{fig:qq_diag}\label{fig1}
\end{figure}

\begin{figure}[t]

\includegraphics{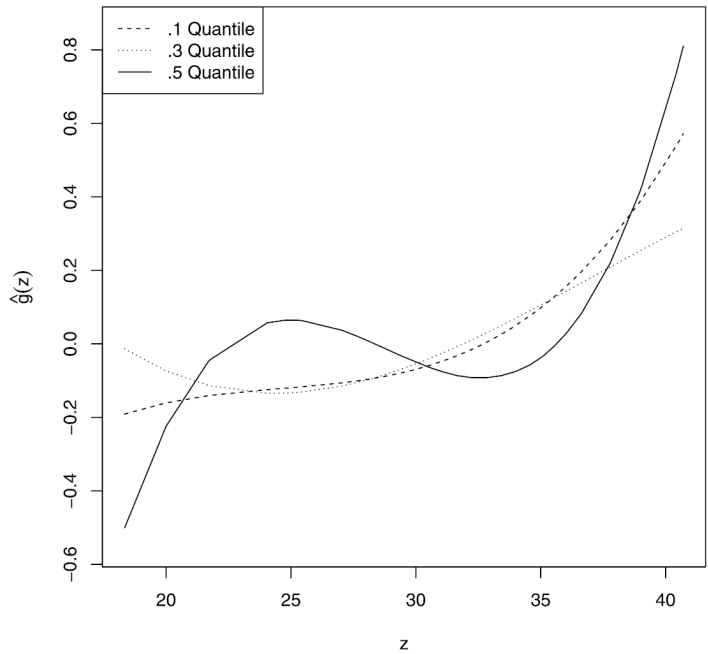}

\caption{Estimated nonlinear effects of mother's age (denoted by $z$)
at three different quantiles.}\label{fig:nl_plots}\label{fig2}
\end{figure}

Model checking in high dimension is challenging. In the following, we
consider a simulation-based diagnostic plot to help visually assess the
overall lack-of-fit for the quantile regression model
[Wei and He (\citeyear{WH06})] to assess the overall lack-of-fit for the quantile regression
model. First, we randomly generate $\tilde{\tau}$ from the uniform
$[0,1]$ distribution. Then we fit the proposed semiparametric quantile
regression model using the SCAD penalty for the quantile $\tilde{\tau
}$. Next, we generate a response variable $\tilde{Y} = \vx'\hat
{\vbeta}(\tilde{\tau}) + \hat{g}(z,\tilde{\tau})$,
where $(\vx,z)$ is randomly sampled from the set of observed
covariates, with $z$ denoting mother's age and
$\vx$ denoting the vector of other covariates. The process is repeated
100 times and produces a sample of 100 simulated birth weights based on
the model. Figure~\ref{fig1} shows the QQ plot comparing the simulated and
observed birth weights. Overall, the QQ plot is close to the 45 degree
line and does not suggest gross lack-of-fit. Figure~\ref{fig2} displays the
estimated nonlinear effects of mother's age $\hat{g}(z)$ at the three
quantiles [standardized to satisfy the constraint $\sum_{i=1}^n \hat
{g}(z_i) = 0$]. At the 0.1 and 0.3 quantiles, the estimated mother's
age effects are similar except for some deviations at the tails of the
mother's age distribution. At these two quantiles, after age 30,
mother's age is observed to have a positive effect. The effect of
mother's age at the median is nonmonotone: the effect is first
increasing (up to age 25), then decreasing (to about age 33), and
increasing again.

We observe that different models are often selected for different
random partitions. Table~\ref{tab6} summarizes the variables selected by Q-SCAD
for $\tau= 0.1$, 0.3 and 0.5 and the frequency these variables are
selected in the 100 random partitions. Probes are listed by their
identification number along with corresponding gene in parentheses. The
SCAD models tend to produce sparser models while the LASSO models
provide slightly better predictive performance.

%

%
\begin{table}[t]
\tabcolsep=0pt
\caption{Frequency of covariates selected at three quantiles among 100 random partitions}\label{table_partition_selections}\label{tab6}
\begin{tabular*}{\tablewidth}{@{\extracolsep{\fill}}@{}lclclc@{}}
\hline
\multicolumn{2}{@{}c}{\textbf{Q-SCAD 0.1}} &
\multicolumn{2}{c}{\textbf{Q-SCAD 0.3}} &
\multicolumn{2}{c@{}}{\textbf{Q-SCAD 0.5}} \\[-6pt]
\multicolumn{2}{@{}c}{\hrulefill} &
\multicolumn{2}{c}{\hrulefill} &
\multicolumn{2}{c@{}}{\hrulefill} \\
 & \textbf{Fre-} &  & \textbf{Fre-} &  & \textbf{Fre-} \\
\textbf{Covariate} & \textbf{quency} & \multicolumn{1}{c}{\textbf{Covariate}} & \textbf{quency} & \multicolumn{1}{c}{\textbf{Covariate}} & \textbf{quency} \\
\hline
Gestational age & 82 & Gestational age & 86 & Gestational age & 69 \\
1,687,073 (SOGA1) & 24 & 1,804,451 (LEO1) & 33 & 2,334,204 (ERCC6L) & 57 \\
& & 1,755,657 (RASIP1) & 27 & 1,732,467 (OR2AG1) & 52 \\
& & 1,658,821 (SAMD1) & 23 & 1,656,361 (LOC201175) & 31 \\
& & 2,059,464 (OR5P2) & 14 & 1,747,184 (PUS7L) & \phantom{0}5 \\
& & 2,148,497 (C20orf107) & \phantom{0}6 & & \\
& & 2,280,960 (DEPDC7) & \phantom{0}3 & &\\
\hline
\end{tabular*}\vspace*{-3pt}
\end{table}

Gestational age is identified to be important with high frequency at
all three quantiles under consideration.
This is not surprising given the known important relationship between
birth weight and
gestational age. Premature birth is often strongly associated with low
birth weight.
The genes selected at the three different quantiles are not
overlapping. This is an
indication of the heterogeneity in the data. The variation in frequency
is likely
due to the relatively small sample size.
However, examining the selected genes does provide some interesting
insights. The gene SOGA1 is a suppressor of glucose, which is
interesting because maternal gestational diabetes is known to have a
significant effect on birth weight [\citet{gestDiab}]. The genes
OR2AG1, OR5P2 and DEPDC7 are all located on chromosome 11,
the chromosome with the most selected genes. Chromosome 11 also
contains PHLDA2, a gene that has been reported to be highly expressed
in mothers that have children with lower birth weight [\citet{phlda2}].

\section{Estimation and variable selection for multiple quantiles}\label{sec6}
Motivated by referees' suggestions, we consider an extension for
simultaneous variable selection at multiple quantiles. Let $\tau_1 <
\tau_2 < \cdots< \tau_M$ be the set of quantiles of interest, where
$M > 0$ is a positive integer. We assume that
%
%
\begin{equation}
\label{multTau} Q_{Y_i \mid\vx_i, \vz_i}(\tau_m) = \vx_i'
\vbeta_0^{(m)} + g_0^{(m)}(
\vz_i),\qquad m = 1,\ldots,M,
\end{equation}

\noindent  where $g_0^{(m)}(\vz_i) = g_{00}^{(m)} + \sum_{j=1}^d
g_{0j}^{(m)}(z_{ij})$, with $g_{00}^{(m)} \in\mathcal{R}$. We assume
that functions $g_{0j}^{(m)}$ satisfy $E
[g_{0j}^{(m)}(z_{ij}) ]=0$ for the purpose of identification. The
nonlinear functions
are allowed to vary with the quantiles. We are interested in the high
dimensional case where
most of the linear covariates have zero coefficients across all $M$
quantiles, for which group selection will help us combine information
across quantiles.

We write $\vbeta_0^{(m)} = (\beta_{01}^{(m)},\beta_{02}^{(m)},
\ldots, \beta_{0p_n}^{(m)} )'$, $m = 1,\ldots,M$.
Let $\bar{\vbeta}^{0j}$ be the \mbox{$M$-}vec\-tor $ (\beta
_{0j}^{(1)},\ldots, \beta_{0j}^{(M)} )'$, $1 \leq j \leq p_n$.
Let $\bar{A} = \{j: \llVert \bar{\vbeta}^{0j}\rrVert \neq0, 1 \leq j
\leq
p_n \}$ be the index set of variables
that are active at least one quantile level of interest, where $\llVert
\cdot\rrVert $ denotes the $L_2$ norm. Let $\bar{q}_n = |\bar{A}|$ be the
cardinality of $A$. Without loss of generality, we assume $\bar{A} = \{
1,\ldots,\bar{q}_n\}$. Let $X_{\bar{A}}$ and $\vx_{\bar
{A}_1}',\ldots,\vx_{\bar{A}_n}$ be defined as before.
By the result of \citet{Schumaker}, there exists $\vxi_0^{(m)} \in
\mathcal{R}^{L_n}$, where \mbox{$L_n = d(k_n+l+1)+1$}, such that
$\mathop{\operatorname{sup}}_{\vz_i} |\bolds{\Pi}(\vz_i)'\vxi
_0^{(m)} -
g_0^{(m)}(\vz_i)| = O(k_n^{-r})$, $m=1,\ldots,M$.

We write the $(Mp_n)$-vector $\vbeta= ({\vbeta^{(1)}}',\ldots
,{\vbeta^{(M)}}' )'$,
where for $k=1,\ldots, M$, $\vbeta^{(k)}=(\beta^{(k)}_1,\ldots
,\beta^{(k)}_{p_n})'$;
and we write the $(ML_n)$-vector $\vxi= ({\vxi^{(1)}}',\ldots
,{\vxi^{(M)}}' )$. Let $\bar{\vbeta}^j$ be the $M$-vector
$ (\beta_j^{(1)},\ldots, \beta_j^{(M)} )'$, $1 \leq j
\leq p_n$. For simultaneous variable selection and estimation, we
estimate $ (\vbeta_0^{(m)},\xi_0^{(m)} )$, $m = 1,\ldots
,M$, by minimizing the following penalized objective function
%
%
\begin{eqnarray}\label{groupPen}
\bar{Q}^P(\vbeta,\vxi) &=& n^{-1} \sum
_{i=1}^n \sum_{m=1}^{M}
\rho_{\tau_m}\bigl(Y_i - \vx_i'
\vbeta^{(m)} - \bolds{\Pi}(\vz_i)'
\vxi^{(m)}\bigr)
\nonumber\\[-8pt]\\[-8pt]\nonumber
&&{} + \sum_{j=1}^{p_n}
p_\lambda\bigl(\bigl\llVert\bar{\vbeta}^{j}\bigr\rrVert
_1\bigr),
\end{eqnarray}
where $p_\lambda(\cdot)$ is a penalty function with tuning parameter
$\lambda$, $\llVert \cdot\rrVert _1$ denotes the $L_1$ norm, which
was used in
\citet{YuanLin} for group penalty; see also \citet
{HuangEtAl}. The penalty function encourages group-wise sparsity and forces
the covariates that have no effect on any of the $M$ quantiles to be
excluded together. Similarly penalty functions have been used in
\citet{zouYuan}, \citet{LiuWu}
for variable selection at multiple quantiles. The above estimator can
be computed similarly as in Section~\ref{sec3.2}.

In the oracle case, the estimator would be obtained by considering the
unpenalized part of (\ref{groupPen}), but with $\vx_i$ replaced by
$\vx_{\bar{A}_i}$. 
That is, we let
%
%
\begin{eqnarray}
\label{mm} && \bigl\{\hat{\vbeta}_1^{(m)},\hat{
\xi}^{(m)}: 1\leq m\leq M \bigr\}
\nonumber\\[-8pt]\\[-8pt]\nonumber
&&\qquad = \mathop{\operatorname{argmin}}_{\vbeta_1^{(m)},\vxi^{(m)},
1\leq m\leq M} n^{-1} \sum
_{i=1}^n \sum_{m=1}^M
\rho_{\tau_{m}}\bigl(Y_i - \vx_{\bar{A}_i}'
\vbeta_1^{(m)} - \bolds{\Pi}(\vz_i)'
\vxi^{(m)}\bigr).
\end{eqnarray}
The oracle estimator for $\vbeta_0^{(m)}$ is $\hat{\vbeta
}^{(m)}= (\hat{\vbeta}_1^{(m)'},\mathbf{0}_{p_n-q_n}' )'$,
and across all quantiles is $\hat{\bar{\vbeta}} = (\hat
{\vbeta}^{(1)},\ldots,\hat{\vbeta}^{(M)} )$ and $\hat{\bar
{\vxi}} = (\hat{\vxi}^{(1)},\ldots,\hat{\vxi}^{(M)}
)$. The oracle estimator for the nonparametric function $g_{0j}^{(m)}$
is $\hat{g}_j^{(m)}(z_{ij}) = \vpi(z_{ij})'\hat{\vxi_j}^{(m)} -
n^{-1}\* \sum_{i=1}^n \pi(z_{ij})'\hat{\vxi}_j^{(m)}$ for $j=1,\ldots
,d$; for $g_{00}^{(m)}$ is $\hat{g}_0^{(m)} = \hat{\xi}_0^{(m)} +
n^{-1} \sum_{i=1}^n \times\break \sum_{j=1}^d \pi(z_{ij})'\hat{\vxi}_j^{(m)}$.
The oracle estimator of $g_0^{(m)}(\vz_i)$ is $\hat{g}^{(m)}(\vz_i)
= \hat{g}_0^{(m)} +\break \sum_{j=1}^d \hat{g}_j^{(m)}(z_{ij})$.
As the next theorem suggests, Theorem \ref{scad_local_min} can be extended to the
multiple quantile case. To save space, we present the regularity
conditions and the technical derivations in the online supplementary
material [\citet{Supp}].

%
\begin{theorem}
\label{scad_local_min_mult}
Assume Conditions \textup{B1--B6} in the online supplementary material [\citet{Supp}] are satisfied.
Let\vspace*{1pt} $\bar{\mathcal{E}}_n(\lambda)$ be the set of local minima of the
penalized objective function $\bar{Q}^P(\vbeta,\gamma)$.
Consider either the SCAD or the MCP penalty\vspace*{1pt} function with tuning
parameter $\lambda$. Let $\hat{\bar{\veta}} \equiv(\hat
{\bar{\vbeta}},\hat{\bar{\vxi}} )$ be the oracle estimator
that solves (\ref{mm}). If 
$\lambda= o (n^{-(1-C_4)/2} )$, $n^{-1/2}\bar{q}_n =
o(\lambda)$, $n^{-1/2}k_n = o(\lambda)$ and $\log(p_n) = o(n\lambda^2)$,
then
\[
P \bigl(\hat{\bar{\veta}} \in\bar{\mathcal{E}}_n(\lambda) \bigr)
\rightarrow1\qquad\mbox{as } n \rightarrow\infty.
\]
\end{theorem}


\subsection*{A numerical example}
To assess the multiple quantile
estimator, we ran 100 simulations using the setting presented in
Section~\ref{sec4} with $\varepsilon_i \sim T_3$, and consider $\tau=
0.5$, 0.7
and 0.9. We compare the variable selection performance of
the multiple-quantile estimator (denoted by Q-group) in this section
with the method that estimates each quantile separately (denoted by
Q-ind). For both approaches, we use the SCAD penalty function. Results
for the MCP penalty are included in the online supplementary material
[\citet{Supp}]. We also report results from the
multiple-quantile oracle estimator (denotes by Q-oracle) which assumes
the knowledge of the underlying model and serves as a benchmark.

%
\begin{table}[b]
\tabcolsep=0pt
\caption{Comparison of group and individual penalty functions for
multiple quantile estimation with $\varepsilon\sim T_3$}\label{tab7}
\begin{tabular*}{\tablewidth}{@{\extracolsep{\fill}}@{}lcd{1.2}d{1.2}d{1.2}c@{}}
\hline
\textbf{Method} & \multicolumn{1}{c}{$\bolds{p}$} & \multicolumn{1}{c}{\textbf{FV}} & \multicolumn{1}{c}{\textbf{TV}} & \multicolumn{1}{c}{\textbf{True}}
& \multicolumn{1}{c@{}}{$\bolds{L_2}$ \textbf{error}}
\\
\hline
Q-group-SCAD & 300 & 1.01 & 4 & 0.49 & 0.14 \\
Q-ind-SCAD & 300 & 0.98 & 4 & 0.45 & 0.17 \\
Q-oracle & 300 & 0 & 4 & 1 & 0.06
\\[3pt]
Q-group-SCAD & 600 & 1.2 & 4 & 0.56 & 0.15 \\
Q-ind-SCAD & 600 & 1.51 & 3.99 & 0.34 & 0.17 \\
Q-oracle & 600 & 0 & 4 & 1 & 0.07 \\
\hline
\end{tabular*}
\end{table}

Table~\ref{tab7} summarizes the simulation results for $n=50$, $p=300$ and 600.
As in \citet{zouYuan}, when evaluating the Q-ind method,
at quantile level $\tau_m$, we define $A_m=\{j: \hat{\beta
}_j^{(m)}\neq0\}$ be
the index set of estimated nonzero coefficients at this quantile level.
Let $\bigcup_{m=1}^MA_m$ be the set of the selected variables using Q-ind.
As the simulations results in Section~\ref{sec4}, we report FV, TV and TRUE.
We also report the error for estimating the linear coefficients ($L_2$ error),
which is defined as the average of $M^{-1} \sum_{m=1}^M (\hat
{\vbeta}^{(m)}-\vbeta_{0}^{(m)} )^2$
over all simulation runs. The results demonstrate that
comparing with Q-ind, the new method Q-group has lower false discovery rate,
higher probability of identifying the true underlying model and smaller
estimation error.

\section{Discussion}\label{sec7}
We considered nonconvex penalized estimation for partially linear
additive quantile regression models
with high dimensional linear covariates. We derive the oracle theory
under mild conditions.
We have focused on estimating a particular quantile of interest and
also considered an extension
to simultaneous variable selection at multiple quantiles.

A problem of important practical interest is how to identify which
covariates should be modeled linearly and which
covariates should be modeled nonlinearly. Usually, we do not have such
prior knowledge in real data analysis.
This is a challenging problem in high dimension. Recently, important
progresses have been made by
\citet{ZCL}; \citet{HWM}; \citet{Lian} for
semiparametric mean
regression models. We plan on addressing this question for high
dimensional semiparametric quantile regression in our future research.

Another relevant problem of practical interest is to estimate the
conditional quantile function itself.
Given $\vx^*$, $\vz^*$, we can estimate $Q_{Y_i \mid\vx^*, \vz
^*}(\tau)$
by ${\vx^*}'\hat{\vbeta}_{1} + \hat{g}(\vz^*)$, where
$\hat{\vbeta}$ and $\hat{g}$ are obtained from penalized quantile regression.
We conjecture that the consistency of estimating the conditional quantile
function can be derived under somewhat weaker conditions
in the current paper, as motivated by the results on \textit{persistency}
for linear mean regression in
high dimension [\citet{Green}]. The details will also be
further investigated
in the future.

\begin{appendix}\label{sec8}\label{appe}
\section*{Appendix}

Throughout the appendix, we use $C$ to denote a positive constant which
does not depend on $n$ and may vary from
line to line. For a vector $\vx$, $\llVert \vx\rrVert $
denotes its Euclidean norm.
For a matrix $A$, $\llVert A\rrVert = \sqrt{\lambda_{\max}(A'A)}$
denotes its
spectral norm. For a function $h(\cdot)$ on $[0,1]$, $\llVert
h\rrVert _{\infty}
= \mathop{\operatorname{\sup}}_{x} |h(x)|$
denotes the uniform norm. Let $I_n$ denote an $n\times n$ identity matrix.

\subsection{Derivation of the results in Section~\texorpdfstring{\protect\ref{sec2}}{2}}\label{sec8.1}

\subsubsection{Notation}\label{sec8.1.1}
To facilitate the proof, we will make use of the theoretically centered
B-spline basis functions
similar to the approach used by \citet{XY}. More specifically,
we consider the B-spline basis functions $b_j(\cdot)$ in\vspace*{2pt} Section~\ref{sec2.1}
and let $B_j(z_{ik}) = b_{j+1}(z_{ik}) - \frac{E
[b_{j+1}(z_{ik}) ]}{E [b_1(z_{ik}) ]}b_1(z_{ik})$
for $j=1,\ldots, k_n+l$. Then $E(B_j(z_{ik}))=0$.
For a given covariate $z_{ik}$,
let
${\mathbf w}(z_{ik})= (B_1(z_{ik}),\ldots,B_{k_n+l}(z_{ik}) )'$
be the vector of basis functions, and $\vW(\vz_i)$ denote the
$J_n$-dimensional vector
$ (k_n^{-1/2},{\mathbf w}(z_{i1})',\ldots,{\mathbf w}(z_{id})' )'$,
where $J_n=d(k_n+l)+1$.

By the result of \citeauthor{Schumaker} [(\citeyear{Schumaker}),
page~227], there exists a
vector $\vgamma_0\in\mathcal{R}^{J_n}$ and a
positive constant
$C_0$, such that $\sup_{t\in[0,1]^{d}}|g_0(\vt)-\vW(\vt)'\vgamma
_0|\leq C_0k_n^{-r}$.
Let
%
%
\begin{equation}
\label{orcObjFun2} (\hat{\bolds{\vc}}_1, \hat{\vgamma} ) = \mathop{
\operatorname{argmin}}_{ (\vc_1, \vgamma)} \frac
{1}{n} \sum
_{i=1}^n\rho_\tau\bigl(Y_i
- \vx_{A_i}'\vc_1 - \vW(\vz
_{i})'\vgamma\bigr).
\end{equation}
We write $\vgamma=(\gamma_0,\vgamma_1',\ldots,\vgamma_d')'$, where
$\gamma_0\in\mathcal{R}$, $\vgamma_j\in\mathcal{R}^{k_n+l}$,
$j=1,\ldots,d$;
and we write $\hat{\vgamma}=(\hat{\gamma}_0,\hat{\vgamma
}_1',\ldots,\hat{\vgamma}_d')'$
the same fashion.
It can be shown that (see the supplemental material) $\hat{\bolds
{c}}_1=\hat{\bolds{\vbeta}}_1$.
So the change of the basis functions for the nonlinear part does not
alter the estimator for the linear part.
Let
$\tilde{g}_j(\vz_i) = w(z_{ij})'\hat{\vgamma}_j$
be the estimator of $g_{0j}$,
$j=1,\ldots,d$. The estimator for $g_{00}$ is
$\tilde{g}_{0}=k_n^{-1/2}\hat{\vgamma}_0$.
The estimator for $g_{0}(\vz_i)$ is $\tilde{g}(\vz_i)=\vW(\vz
_i)'\hat
{\vgamma}=\tilde{g}_0+
\sum_{j=1}^{d}\tilde{g}_j(\vz_i)$. It can be derived that (see the
supplemental material)
$\hat{g}_j(\vz_i)=\tilde{g}_j(\vz_i)-n^{-1}\sum_{i=1}^n\tilde
{g}_j(\vz_i)$
and $\hat{g}_0=\tilde{g}_{0}+n^{-1}\sum_{i=1}^n\sum_{j=1}^d\tilde
{g}_j(\vz_i)$.
Hence, $\hat{g}=\hat{g}_0+\sum_{j=1}^d\hat{g}_j=\tilde{g}$.
Later, we will show
$n^{-1} \sum_{i=1}^n (\tilde{g}(\vz_i) - g_{0}(\vz_i)
)^2 =
O_p (n^{-1}(q_n+dJ_n) )$.


Throughout the proof, we will also use the following notation:
\begin{eqnarray*}
\psi_\tau(\varepsilon_i) &=& \tau-I(\varepsilon_i < 0),
\\
W&=& \bigl(\vW(\vz_1),\ldots,\vW(\vz_n)
\bigr)' \in\mathds{R}^{n \times
J_n},
\\
P &=& W\bigl(W'B_nW\bigr)^{-1}W'B_n
\in\mathds{R}^{n \times n},
\\
X^* &=&\bigl(\vx_1^*,\ldots,\vx_n^*
\bigr)'=(I_n-P)X_A \in\mathds{R}^{n
\times q_n},
\\
W_B^2 &=& W'B_nW \in
\mathds{R}^{J_n \times J_n},
\\
\vtheta_1 &=& \sqrt{n} (\vc_1-
\vbeta_{10} )\in\mathds{R}^{q_n},
\\
\vtheta_2 &=& W_B (\vgamma-\vgamma_0 ) +
W_B^{-1}W'B_nX_A(
\vc_1-\vbeta_{10}) \in\mathds{R}^{J_n},
\\
\tilde{\vx}_i &=& n^{-1/2}\vx_i^* \in
\mathds{R}^{q_n},
\\
\tilde{\vW}(\vz_i)&=& W_B^{-1}\vW(
\vz_i)\in\mathds{R}^{J_n},
\\
\tilde{\vs}_i &=& \bigl(\tilde{\vx}_i',
\tilde{\vW}(\vz_i) \bigr)' \in\mathds{R}^{q_n+J_n},
\\
u_{ni}&=&\vW(\vz_i)'\vgamma_{0}
- g_{0}(\vz_i).
\end{eqnarray*}

Notice that
\[
n^{-1}\sum_{i=1}^n
\rho_\tau\bigl(Y_i - \vx_{A_i}'
\vc_1 - \vW(\vz_i)'\vgamma\bigr) =
n^{-1}\sum_{i=1}^n
\rho_\tau\bigl(\varepsilon_i - \tilde{\vx}_i'
\vtheta_1 - \tilde{\vW}(\vz_i)'
\vtheta_2 - u_{ni}\bigr).
\]
%
Define the minimizers under the transformation as
\[
(\hat{\vtheta}_1,\hat{\vtheta}_2 ) = \mathop{
\operatorname{arg}\operatorname{min}}_{(\vtheta_1,\vtheta_2)} n^{-1}
\sum
_{i=1}^n\rho_\tau\bigl(
\varepsilon_i - \tilde{\vx}_i'\vtheta
_1 - \tilde{\vW}(\vz_i)'
\vtheta_2 - u_{ni}\bigr).
\]
Let $a_n$ be a sequence of positive numbers and define
\begin{eqnarray*}
Q_i(a_n) &\equiv& Q_i(a_n
\vtheta_1,a_n\vtheta_2) =
\rho_\tau\bigl(\varepsilon_i - a_n\tilde{
\vx}_i'\vtheta_1 - a_n\tilde{
\vW}(\vz_i)'\vtheta_2 - u_{ni}
\bigr),
\\
E_s[Q_i] &=& E [Q_i
\mid\vx_i,\vz_i ].
\end{eqnarray*}

Let $\vtheta=(\vtheta_1', \vtheta_2')'$. Define
%
%
\begin{eqnarray}\label{d_i_defined}
D_i(\vtheta,a_n) &=&
Q_i(a_n) - Q_{i}(0) - E_s
\bigl[Q_i(a_n)-Q_i(0) \bigr]
\nonumber\\[-8pt]\\[-8pt]\nonumber
&&{} + a_n \bigl(\tilde{\vx}_i'
\vtheta_1 + \tilde{\vW}(\vz_i)'\vtheta
_2 \bigr)\psi_\tau(\varepsilon_i).
\end{eqnarray}
%
Noting that $\rho_\tau(u) = \frac{1}{2}|u| + (\tau- \frac
{1}{2} )u$, we have
%
%
\begin{eqnarray}
\label{q_i_equality} Q_i(a_n) - Q_i(0) &=&
\tfrac{1}{2} \bigl[ \bigl\llvert\varepsilon_i - a_n
\tilde{\vx}_i'\vtheta_1 - a_n
\tilde{\vW}(\vz_i)'\vtheta_2 -
u_{ni} \bigr\rrvert- \llvert\varepsilon_i -
u_{ni} \rrvert\bigr]
\nonumber\\[-8pt]\\[-8pt]\nonumber
&&{}- a_n \bigl(\tau- \tfrac{1}{2} \bigr) \bigl(\tilde{
\vx}_i'\vtheta_1 + \tilde{\vW}(
\vz_i)'\vtheta_2 \bigr).
\end{eqnarray}
Define
\[
Q_i^*(a_n) = \tfrac{1}{2} \bigl[ \bigl\llvert
\varepsilon_i - \tilde{\vx}_i'
\vtheta_1a_n - \tilde{\vW}(\vz_i)'
\vtheta_2a_n - u_{ni} \bigr\rrvert- \llvert
\varepsilon_i - u_{ni} \rrvert\bigr].
\]
Then by combining (\ref{d_i_defined}) and (\ref{q_i_equality}),
%
%
\begin{equation}
\label{d_i_defined_2} D_i(\vtheta,a_n) = Q_i^*(a_n)
- E_s\bigl[Q_i^*(a_n)\bigr] +
a_n \bigl(\tilde{\vx}_i'
\vtheta_1 + \tilde{\vW}(\vz_i)'
\vtheta_2 \bigr)\psi_\tau(\varepsilon_i).
\end{equation}
%

\subsubsection{Some technical lemmas}\label{sec8.1.2}

The proofs of Lemmas \ref{spline}--\ref{lem_g_hat_rate} below are
given in the supplemental material [\citet{Supp}].

%
\begin{lemma}\label{spline}
We have the following properties for the spline basis vector:
\begin{longlist}[(3)]
\item[(1)] $E(\llVert \vW(\vz_i)\rrVert )\leq b_1$, $\forall i$, for some
positive constant
$b_1$ for
all $n$ sufficiently large.

\item[(2)] There exists positive constant $b_2$ and $b_2^*$ such that for
all $n$ sufficiently large
$E(\lambda_{\min}(\vW(\vz_i)\vW(\vz_i)^T))\geq b_2k_n^{-1}$ and
$E(\lambda_{\max}(\vW(\vz_i)\vW(\vz_i)^T))\leq b_2^*k_n^{-1}$.

\item[(3)] $E(\llVert W_B^{-1}\rrVert )\geq b_3\sqrt{k_nn^{-1}}$,
for some positive constant $b_3$, for all $n$ sufficiently large.

\item[(4)]
$\mathop{\max}_{i} \llVert \tilde{\vW}(\vz_i)\rrVert =O_p(\sqrt{\frac
{k_n}{n}})$.

\item[(5)] $\sum_{i=1}^nf_i(0)\tilde{\vx}_i\tilde{\vW}(\vz_i)'=\vnull$.
\end{longlist}
\end{lemma}

%
\begin{lemma}
\label{x_star_big_q}
If Conditions \ref{cond_f}--\ref{cond_sigma_large_p} are satisfied,
then:
\begin{longlist}[(3)]
\item[(1)] There exists a positive constant $C$ such that
$\lambda_{\max} (n^{-1} {X^*}'X^* ) \leq C$, with
probability one.

\item[(2)] $n^{-1/2}X^* = n^{-1/2}\Delta_n + o_p(1)$.
Furthermore, $n^{-1}{X^*}'B_nX^* = K_n + o_p(1)$, where $B_n$ and $K_n$
are defined as in Theorem \ref{large_q_clt}.
\end{longlist}
\end{lemma}

%
\begin{lemma}
\label{lem_g_hat_rate}
If Conditions \ref{cond_f}--\ref{cond_sigma_large_p} hold, then
$n^{-1}\sum_{i=1}^n(\tilde{g}(\vz_i)-g_0(\vz_i))^2 = O_p
(d_n/n )$.
\end{lemma}
%


%
\begin{lemma}
\label{theta_tilde_rate}
Assume Conditions \ref{cond_f}--\ref{cond_sigma_large_p} hold.
Let $\tilde{\vtheta}_1 = \sqrt{n} ({X^*}'B_nX^*
)^{-1}\* {X^*}'\psi_\tau(\varepsilon)$,
where $\psi_\tau(\varepsilon)=(\psi_\tau(\varepsilon_1),\ldots,\psi
_\tau(\varepsilon_n))'$,
then:
\begin{longlist}[(3)]
\item[(1)] $\llVert \tilde{\vtheta}_1\rrVert = O_p (\sqrt{q_n} )$.

\item[(2)]
$A_n\Sigma_n^{-1/2}\tilde{\vtheta}_1 \convInD N(0, G)$,
where $A_n$, $\Sigma_n$ and $G$ are defined in Theorem \ref{large_q_clt}.
\end{longlist}
\end{lemma}

\begin{pf}
(1) The result follows from the observation that, by Lemma \ref{x_star_big_q},
\[
\tilde{\vtheta}_1 = \bigl(K_n+o_p(1)
\bigr)^{-1} \bigl[n^{-1/2}\Delta_n'
\psi_\tau(\varepsilon)+n^{-1/2}(H-PX_A)\psi
_\tau(\varepsilon) \bigr],
\]
and $n^{-1/2}\llVert H-PX_A\rrVert =o_{p}(1)$.

(2)
%
\begin{eqnarray*}
A_n\Sigma_n^{-1/2}\tilde{\vtheta}_1
&=& A_n\Sigma_n^{-1/2}K_n^{-1}
\bigl[n^{-1/2}\Delta_n'\psi_\tau(
\varepsilon) \bigr]\bigl(1+o_p(1)\bigr)
\\
&&{} +A_n\Sigma_n^{-1/2}K_n^{-1}
\bigl[n^{-1/2}(H-PX_A) \bigr]\psi_\tau(\varepsilon)
\bigl(1+o_p(1)\bigr),
\end{eqnarray*}
where the second term is $o_p(1)$ because $n^{-1/2}\llVert
H-PX_A\rrVert =o(1)$.
We write $A_n\Sigma_n^{-1/2}K_n^{-1} [n^{-1/2}\Delta_n'\psi_\tau
(\varepsilon) ]
=\sum_{i=1}^nD_{ni}$, where
\[
D_{ni}=n^{-1/2}A_n\Sigma_n^{-1/2}K_n^{-1}\vdelta_i \psi_\tau
(\varepsilon_i).
\]
To verify asymptotic normality, we first note that $E(D_{ni})=0$ and
\begin{eqnarray*}
\sum_{i=1}^nE\bigl(D_{ni}D_{ni}'
\bigr) &=&A_n\Sigma_n^{-1/2}K_n^{-1}S_nK_n^{-1}
\Sigma_n^{-1/2}A_n'=
A_nA_n'\rightarrow G.
\end{eqnarray*}

The proof is complete by
checking the Lindeberg--Feller condition. For any
$\varepsilon> 0$ and using Conditions \ref{cond_f}, \ref{highd_cond_x}
and \ref{cond_sigma_large_p}
\begin{eqnarray*}
&&\sum_{i=1}^nE \bigl[\llVert
D_{ni}\rrVert^2 I\bigl(\llVert D_{ni}\rrVert>
\varepsilon\bigr) \bigr]
\\
&&\qquad \leq\varepsilon^{-2} \sum
_{i=1}^nE\llVert D_{ni}\rrVert
^4
\\
&&\qquad\leq(n\varepsilon)^{-2} \sum_{i=1}^nE
\bigl(\psi^4_\tau(\varepsilon_i) \bigl(
\vdelta_i'K_n^{-1}
\Sigma_n^{-1/2} A_n'A_n
\Sigma_n^{-1/2}K_n^{-1}
\vdelta_i \bigr)^2 \bigr)
\\
&&\qquad\leq C n^{-2}\varepsilon^{-2}\sum
_{i=1}^nE\bigl(\llVert\vdelta_i\rrVert
^4\bigr)= O_p\bigl(q_n^2/n
\bigr)= o_p(1),
\end{eqnarray*}
where the last inequality follows by observing that
$\lambda_{\max}(A_n'A_n)=\lambda_{\max}(A_nA_n')\rightarrow c$ for some finite
positive constant $c$.
\end{pf}

%
%
\begin{lemma}
\label{lem_theta_tilde_hat}
If Conditions \ref{cond_f}--\ref{cond_sigma_large_p} hold, then
\[
\llVert\hat{\vtheta}_1 - \tilde{\vtheta}_1\rrVert=
o_p(1).
\]
\end{lemma}

\begin{pf}
Proof provided in online supplementary material [\citet{Supp}].
\end{pf}

\subsubsection{Proof of Theorems \texorpdfstring{\protect\ref{large_q_oracle}}{2.1},
\texorpdfstring{\protect\ref{large_q_clt}}{2.2} and Corollary \texorpdfstring{\protect\ref{fixed_q_clt}}{1}}\label
{sec8.1.3}

By the observation $\hat{g}=\tilde{g}$, Lemma \ref{lem_g_hat_rate}
implies the second result of Theorem~\ref{large_q_oracle}. The first result
of Theorem~\ref{large_q_oracle} follows by observing
$\hat{\vc}_1=\hat{\vbeta}_1$ and Lemmas \ref{theta_tilde_rate} and
\ref{lem_theta_tilde_hat}.
The proof of Theorem \ref{large_q_clt} follows from Lemmas \ref
{theta_tilde_rate} and \ref{lem_theta_tilde_hat}.
Set $A_n = I_{q}$, then the proof of Corollary \ref{fixed_q_clt}
follows from the fact that $q$ being constant and Theorems \ref
{large_q_oracle} and \ref{large_q_clt}.


\subsection{Derivation of the results in Section~\texorpdfstring{\protect\ref{sec3.3}}{3.3}}\label{sec8.2}

%
\begin{lemma}
\label{lem_diff_convex}
Consider the function $k(\eta)-l(\eta)$ where both $k$ and $l$ are
convex with subdifferential functions $\partial k(\eta)$ and $\partial
l(\eta)$.
Let $\eta^*$ be a point that has neighborhood $U$ such that $\partial
l(\eta) \cap\partial k(\eta^*) \neq\varnothing, \forall \eta\in
U \cap\operatorname{dom}(k)$.
Then $\eta^*$ is a local minimizer of $k(\eta) - l(\eta)$.
\end{lemma}

\begin{pf}
The proof is available in \citet{TaoAn97}.
\end{pf}


\subsubsection{Proof of Lemma \texorpdfstring{\protect\ref{lem_sub_diff_q}}{1}}\label{sec8.2.1}
\mbox{}

\begin{pf*}{Proof of (\ref{lem_sub_diff_q_1})}
By convex optimization theory $\mathbf{0} \in\partial\sum_{i=1}^n\rho
_\tau
(Y_i - \vx_i'\vbeta- \bolds{\Pi}(\vz_i)'\vxi)$.
Thus, there exists
$a_j^*$ as described in the lemma such that with the choice
$a_j=a_j^*$, we have $s_j(\hat{\vbeta},\hat{\vxi})=0$ for
$j=1,\ldots,q_n$ or $j=p_n+1,\ldots,p_n+J_n$.
\end{pf*}

\begin{pf*}{Proof of (\ref{lem_sub_diff_q_2})}
It is sufficient to show
$
P (\llvert \hat{\beta}_j\rrvert \geq(a+1/2)\lambda$,
for $j= 1,\ldots,q_n ) \rightarrow1
$
as $n,p\rightarrow\infty$.
Note that
%
%
\begin{equation}
\mathop{\min}_{1\leq j \leq q_n} \llvert\hat{\beta}_j\rrvert\geq
\mathop{\min}_{1 \leq j \leq q_n} \llvert\beta_{0j}\rrvert- \mathop{
\max}_{1
\leq j \leq q_n} \llvert\hat{\beta}_j - \beta_{0j}
\rrvert.
\end{equation}
By Condition \ref{cond_small_sig}, $\mathop{\min}_{1 \leq j \leq
q_n} \llvert \beta_{0j}\rrvert \geq C_5 n^{-(1-C_4)/2}$.
By Theorem \ref{large_q_oracle} and Conditions \ref
{cond_sigma_large_p} and \ref{cond_small_sig},
$\mathop{\max}_{1 \leq j \leq q_n} \llvert \hat{\beta}_j - \beta
_{0j}\rrvert = O_p (\sqrt{\frac{q_n}{n}} ) = o_p
(n^{-(1-C_4)/2} )$.
(\ref{lem_sub_diff_q_1}) holds by noting $\lambda= o
(n^{-(1-C_4)/2} )$.
\end{pf*}

\begin{pf*}{Proof of (\ref{lem_sub_diff_q_3})}
Proof provided in the online supplementary material [\citet{Supp}].
\end{pf*}

\subsubsection{Proof of Theorem \texorpdfstring{\protect\ref{scad_local_min}}{3.1}}\label{sec8.2.2}

Recall that for $\kappa_{j} \in\partial k(\vbeta,\vxi)$ 
%
\begin{eqnarray*}
\kappa_j &=& s_j(\vbeta,\vxi) + \lambda l_j
\qquad\mbox{for } 1 \leq j \leq p_n,
\\
\kappa_{j} &=& s_{j}(\vbeta,\vxi) \qquad\mbox{for }
p_n+1 \leq j \leq p_n+J_n.
\end{eqnarray*}

Define the set
\begin{eqnarray*}
\mathcal{G} &=& \bigl\{ \vkappa= (\kappa_1, \kappa_2,
\ldots,\kappa_{p_n+J_n})': \kappa_j = \lambda
\operatorname{sgn}(\hat{\beta}_j), j=1,\ldots,q_n;
\\
&&{} \kappa_j = s_j(\hat{\vbeta}, \hat{\vxi}) + \lambda
l_j, j=q_n+1,\ldots,p_n;
\\
&&{} \kappa_j = 0, j=p_n+1,\ldots,p_n+J_n,
\bigr\},
\end{eqnarray*}
where $l_j$ ranges over $[-1,1]$ for $j= q_n+1,\ldots,p_n$.
By Lemma \ref{lem_sub_diff_q}, we have $P(\mathcal{G} \subset\partial k(\hat{\vbeta
},\hat{\vxi}))\rightarrow1$.

Consider any $(\vbeta', \vxi')'$ in a ball with the center $
(\hat{\vbeta}',\hat{\vxi}' )'$ and radius
$\lambda/ 2$.
By Lemma \ref{lem_diff_convex}, to prove the theorem it is sufficient
to show that there exists $\vkappa^* = (\kappa_1^*,\ldots
,\kappa^*_{p_n+J_n} )' \in\mathcal{G}$
such that
%
%
\begin{eqnarray}
P \biggl( \kappa_j^* = \frac{\partial l(\vbeta, \vxi)}{\partial
\beta_j}, j=1,
\ldots,p_n \biggr) &\rightarrow&1; \label{xi_star_j}
\\
P \biggl( \kappa_{p_n+j}^* = \frac{\partial l(\vbeta, \vxi
)}{\partial\xi_j}, j=1,
\ldots,J_n \biggr) &\rightarrow&1. \label
{xi_star_j2}
\end{eqnarray}

Since $\frac{\partial l(\vbeta, \vxi)}{\partial\xi_j} = 0$ for
$j=1,\ldots,J_n$,
(\ref{xi_star_j2}) is satisfied by Lemma \ref{lem_sub_diff_q}.
We outline how $\kappa_j^*$ can be selected to satisfy (\ref
{xi_star_j}). 
%
\begin{longlist}[2.]
%
\item[1.] For $1 \leq j \leq q_n$,
we have $\kappa^*_j=\lambda\operatorname{sgn}(\hat{\beta}_j)$ for $\beta
_j \neq0$. 
For either SCAD or MCP penalty function, $\frac{\partial l(\vbeta,
\vxi)}{\partial\beta_j} = \lambda\operatorname{sgn}(\beta_j)$ for
$\llvert \beta_j\rrvert > a\lambda$. By Lemma \ref{lem_sub_diff_q},
we have
\begin{eqnarray*}
\min_{1\leq j \leq q_n}\llvert\beta_j\rrvert&\geq& \min
_{1 \leq j \leq
q_n}\llvert\hat{\beta}_j\rrvert- \max
_{1\leq j \leq q_n}\llvert\hat{\beta}_j - \beta_j
\rrvert
\geq (a+1/2)\lambda- \lambda/2 = a\lambda,
\end{eqnarray*}
with probability approaching one.
Thus, $P(\frac{\partial l(\vbeta, \vxi)}{\partial\beta_j} =
\lambda\operatorname{sgn}(\beta_j))\rightarrow1$.
For any $1 \leq j \leq q_n$, $\llVert \hat{\beta}_j-\beta_{0j}\rrVert =
O_p (n^{-1/2}q_n^{1/2} ) = o(\lambda)$. Therefore, for
sufficiently large $n$, $\hat{\beta}_j$ and $\beta_j$ have the same
sign. This implies
$P(\frac{\partial l(\vbeta, \vxi)}{\partial\beta_j} = \kappa_j^*,
1 \leq j \leq q_n)\rightarrow1$ as $n\rightarrow\infty$.
\item[2.]
For $j=q_n+1,\ldots,p_n$,
$\hat{\beta}_j=0$ by the definition of the oracle estimator and
$\kappa_j = \lambda l_j$ with $l_j \in[-1,1]$. Therefore,
\[
\llvert\beta_j\rrvert\leq\llvert\hat{\beta}_j\rrvert
+ \llvert\hat{\beta}_j - \beta_j\rrvert< \lambda/2.
\]
For $\llvert \beta_j\rrvert < \lambda$, $\frac{\partial l(\vbeta, \vxi
)}{\partial\beta_j} = 0$ for the SCAD penalty and $\frac{\partial
l(\vbeta, \vxi)}{\partial\beta_j} = \beta_j / a$ for MCP,
$j=q_n+1,\ldots,p_n$.
Note that for both penalty functions, we have $ |\frac{l(\vbeta,
\vxi)}{\partial\beta_j} | \leq\lambda$,
$j= q_n+1,\ldots,p_n$. By Lemma \ref{lem_sub_diff_q}, $|s_j(\hat
{\beta}_j)| \leq\lambda/2$ with probability approaching one for $j =
q_n+1,\ldots, p_n$. Therefore, for both penalty functions, there
exists $l_j^* \in[-1,1]$ such that $P(s_j(\hat{\vbeta},\hat{\vxi})
+ \lambda l_j^* = \frac{\partial l(\vbeta, \vxi)}{\beta_j}, j =
q_n+1,\ldots, p_n) \rightarrow1$. Define $\kappa_j^* = s_j(\hat
{\vbeta},\hat{\vxi}) + \lambda l_j^*$.
Then $P(\frac{\partial l(\vbeta, \vxi)}{\partial\beta_j} = \kappa
_j^*, q_n+1 \leq j \leq p_n)\rightarrow1$ as $n\rightarrow\infty$.
\end{longlist}
This completes the proof.
\end{appendix}
\section*{Acknowledgments}
We thank the Editor, the Associate Editor and the anonymous referees
for their careful reading
and constructive comments which have helped us to significantly improve
the paper.

\begin{supplement}[id=suppA]
\stitle{Supplemental Material to ``Partially linear additive quantile regression in~ultra-high dimension''}
\slink[doi]{10.1214/15-AOS1367SUPP} 
\sdatatype{.pdf}
\sfilename{aos1367\_supp.pdf}
\sdescription{We provide technical details for some of the proofs and additional simulation results.}
\end{supplement}

%

\printaddresses
\end{document}